\documentclass[11pt]{article}

\usepackage{amssymb,amsthm,amsmath}
\usepackage{graphicx}
\usepackage{psfrag}
\usepackage{hyperref}
\usepackage{fullpage}

\newtheorem{thm}{Theorem}[subsection]
\newtheorem{cor}[thm]{Corollary}
\newtheorem{lem}[thm]{Lemma}

\newtheorem{claim}{Claim}[thm]

\newtheorem{case}{Case}
\newtheorem{claim*}{Claim}
\newtheorem{rem}[thm]{Remark}
\newtheorem{subcase}{Subcase}[case]
\newtheorem{assumpt}[thm]{Assumption}

\newcommand{\bpf}{
\begin{proof}
}
\newcommand{\epf}{
\end{proof}
}

\def\co{\colon\thinspace}

\numberwithin{equation}{section}

\newcommand{\ve}{\varepsilon}

\newcommand{\ovx}{\overline{x}}
\newcommand{\ovy}{\overline{y}}
\newcommand{\ova}{\overline{a}}
\newcommand{\ovc}{\overline{c}}
\newcommand{\ovd}{\overline{d}}
\newcommand{\ove}{\overline{e}}

\newcommand{\ben}{\begin{enumerate}}
\newcommand{\een}{\end{enumerate}}
\newcommand{\bit}{\begin{itemize}}
\newcommand{\eit}{\end{itemize}}
\newcommand{\ita}{\item[(a)]}
\newcommand{\itb}{\item[(b)]}
\newcommand{\itc}{\item[(c)]}
\newcommand{\itd}{\item[(d)]}

\newcommand{\wh}{\widehat}
\newcommand{\wt}{\widetilde}
\newcommand{\mc}{\mathcal}
\newcommand{\intr}{\text{int}\,}
\newcommand{\clo}{\text{cl}\,}

\newcommand{\Set}[2]{\left\{ #1 \ \left| \ #2\right.\right\}}

\newcommand{\grp}[1]{\langle #1\rangle}

\newcommand{\mZ}{\mathbb{Z}}

\newcommand{\Figw}[4]{
\begin{center}
\includegraphics[width=#1]{#2}
\caption{ #3 \label{#4} }
\end{center}}

\begin{document}

\title{Hyperbolic $(1,2)$-knots in $S^3$ with crosscap number
two\\ and tunnel number one}%

\author{Enrique Ram\'{\i}rez-Losada\footnote{
Centro de Investigaci\'on en Matem\'aticas, A.C. Guanajuato,
Gto. 36240, M\'exico; \texttt{kikis@cimat.mx}} \ and  Luis G.
Valdez-S\'anchez \footnote{Corresponding Author: Department of
Mathematical Sciences, University of Texas at El Paso, El Paso,
TX 79968,
USA; \texttt{valdez@math.utep.edu}}}%

\date{}

\maketitle

\begin{abstract}
A knot in $S^3$ is said to have crosscap number two if it
bounds a once-punctured Klein bottle but not a Moebius band. In
this paper we give a method of constructing crosscap number two
hyperbolic $(1,2)$-knots with tunnel number one which are
neither 2-bridge nor (1,1)-knots. An explicit infinite family
of such knots is discussed in detail.
\end{abstract}


\section{Introduction}\label{introduction}

Let $K$ be a knot in the 3-sphere $S^3$ with exterior
$X_K=S^3\setminus\intr N(K)$ (where $N(\cdot)$ denotes regular
neighborhood). For $r$ a slope in $\partial X_K$, we let
$K(r)=X_K\cup_{\partial}S^1\times D^2$ denote the manifold
obtained by performing surgery on $K$ along the slope $r$, so
that $r$ bounds a disk in $S^1\times D^2$. A {\it Seifert Klein
bottle} for $K$ is a once-punctured Klein bottle $P$ properly
embedded in $X_K$ which has integral boundary slope; such a
surface $P$ is {\it unknotted} if $\clo(X_K\setminus N(P))$ is
a genus two handlebody. We say that $K$ has {\it crosscap
number two} if $K$ has a Seifert Klein bottle and its exterior
contains no properly embedded Moebius band. The knot $K$ has
{\it tunnel number one} if there is a properly embedded arc
$\tau$ in $X_K$ such that $\clo(X_K\setminus N(\tau))$ is a
(genus two) handlebody. We also say that the knot $K$ {\it
admits a $(g,n)$ decomposition,} or that $K$ is a {\it
$(g,n)$-knot}, if there is a Heegaard splitting surface $S$ in
$S^3$ of genus $g$ which intersects $K$ transversely and bounds
handlebodies $H, H'$, such that both $K\cap H\subset H$ and
$K\cap H'\subset H'$ are trivial $n$-string arc systems.
Finally, we will use the notation $S^2(a,b,c)$ to represent any
small Seifert fibered space over a 2-sphere with three singular
fibers of indices $a,b,c$.

Any (1,1)-knot has tunnel number one; the converse, however,
does not hold in general. It is therefore remarkable that for
genus one hyperbolic knots the properties of having tunnel
number one, admitting a (1,1) decomposition, or being 2-bridge
are all mutually equivalent; this is the content of the
Goda-Teragaito conjecture, which is the main result of
\cite{scharlemann4}. Since a genus one knot bounds a once
punctured torus, which is the orientable homotopy equivalent of
a once punctured Klein bottle, one might expect crosscap number
two hyperbolic knots to exhibit similar behavior, ie with
having tunnel number one, admitting a (1,1) decomposition, and
being 2-bridge are all equivalent conditions. That this is not
the case follows from
\cite[Theorem 1.1]{valdez8}, which shows that crosscap number two
hyperbolic (1,1)-knots are in general not 2-bridge.

One may still ask if, for crosscap number two hyperbolic knots,
having tunnel number one is equivalent to admitting a (1,1)
decomposition. In this paper we show this is not the case by
constructing explicit examples of crosscap number two
hyperbolic knots that have tunnel number one but do not admit
$(1,1)$ decompositions; moreover, all such knots admit a
$(1,2)$ decomposition. The construction of such examples is
based on the classification of the nontrivial crosscap number
two (1,1)-knots in $S^3$ given in \cite{valdez8}. Explicit
examples of general tunnel number one knots without $(1,1)$
decompositions were first given by Morimoto- Sakuma-Yokota
\cite{sakuma2}, and more recently by Eudave-Mu\~noz
\cite{eudave6} (see also \cite{gordon}).

\begin{figure}
\psfrag{D0}{$D_0$}\psfrag{D1}{$D_1$}
\psfrag{H}{$H$}\psfrag{x}{$x$}
\psfrag{y}{$y$}
\Figw{4.5in}{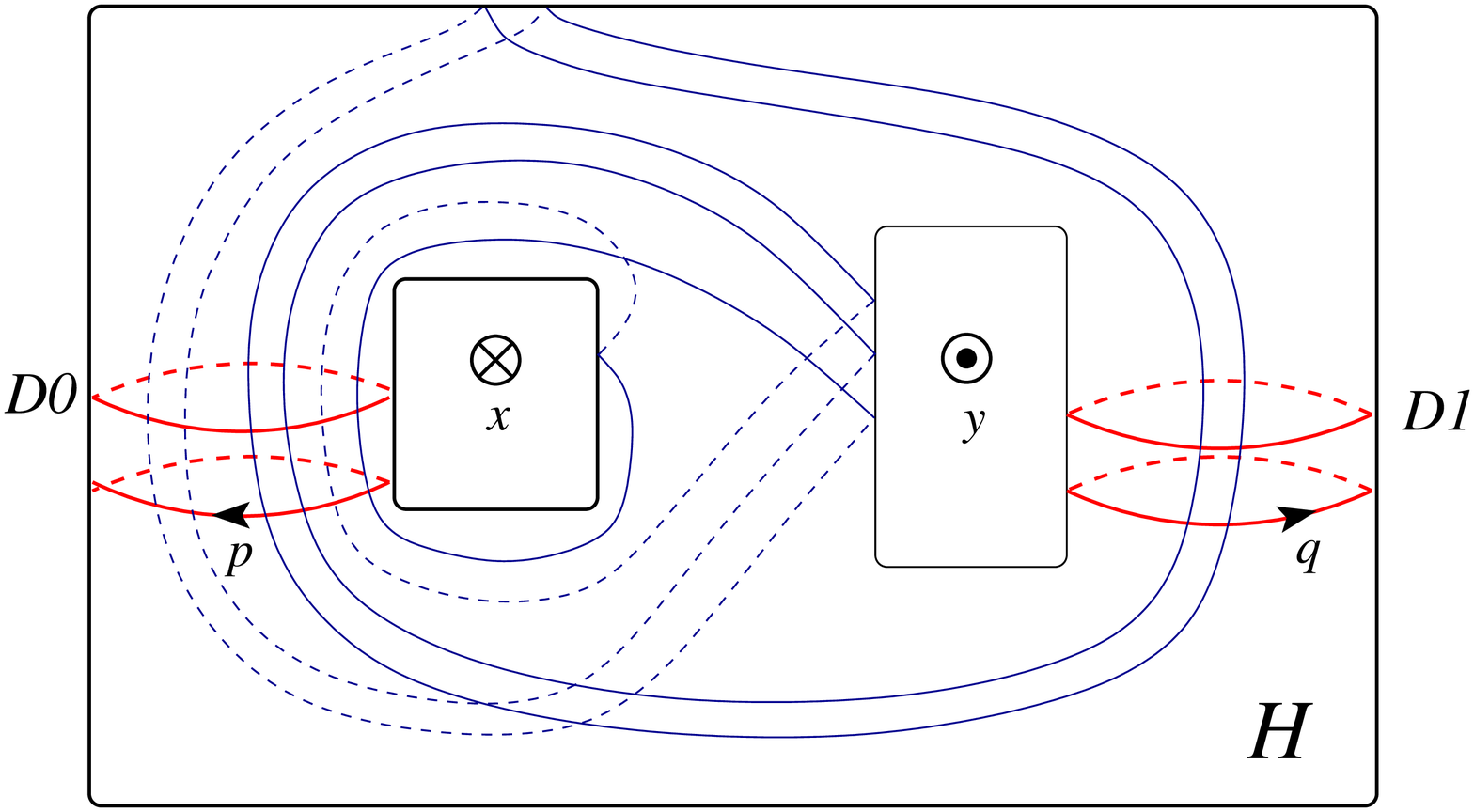}{The knot
$K(0,0)\subset\partial H$.}{fig11b}

\vspace{1.5cm}

\psfrag{D}{$D_0$}\psfrag{p}{$p$}\psfrag{q}{$q=+1$}
\psfrag{T0}{$S_0$}\psfrag{T1}{$S_1$}
\Figw{4in}{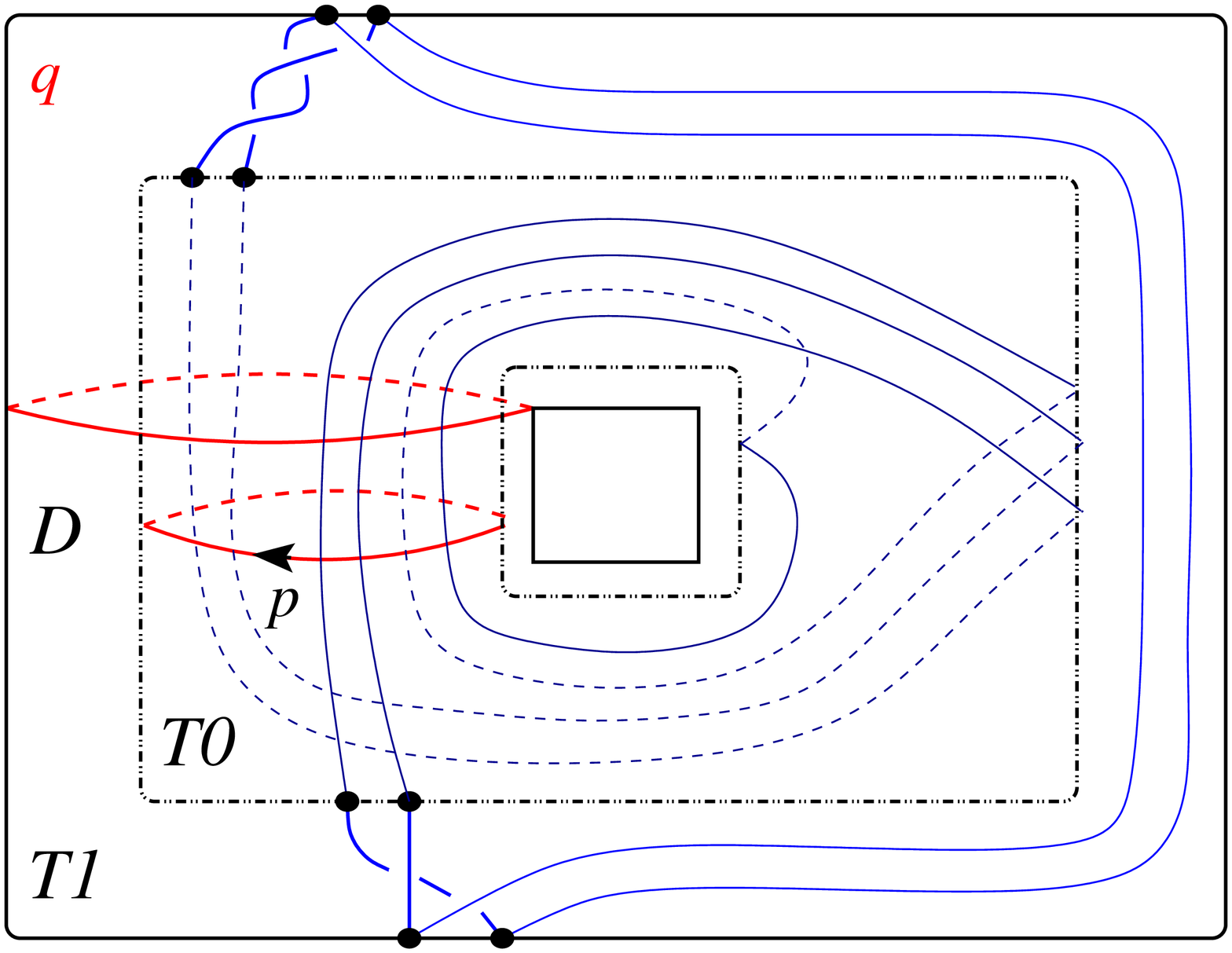}{A $(1,2)$ decomposition of the knot
$K(0,1)$.}{fig11e}
\end{figure}

Our family of examples is constructed starting from the trivial
knot $K(0,0)$ shown in Fig.~\ref{fig11b}, which lies in the
boundary $\partial H$ of an unknotted genus two handlebody $H$
standardly embedded in $S^3$. We obtain a two-parameter
infinite family of knots $K(p,q)\subset\partial H$, for any
integers $p,q$, by Dehn-twisting $K(0,0)$ $p$-times along
$\partial D_0$ and $q$-times along $\partial D_1$, where
$D_0,D_1$ is the complete meridian system for $H$ shown in the
figure; the Dehn-twists are carried out by cutting $\partial H$
along $\partial D_i$, obtaining a 4-punctured 2-sphere $F_0$,
and then twisting only the two `bottom' boundary circles of
$F_0$ indicated by the arrows of the figure the required number
of times (in the direction of the arrows for $p,q>0$, see eg
the knot $K(-1,1)$ in Fig.~\ref{a11b2}).

It is not hard to see that each knot $K(p,q)$ is a
$(1,2)$-knot. Fig.~\ref{fig11e} shows a pair of parallel
unknotted tori $S_0,S_1$ in $S^3$, with the region between
$S_0$ and $S_1$ a product of the form $S_0\times I$; in the
figure, the knot $K(0,1)$ is represented as the union of 4
pairs of arcs: one pair on each $S_0,S_1$ and two pairs of arcs
in $S_0\times I$, with each arc in the latter pairs
intersecting each level torus $S_0\times\{t\}$, $t\in I$,
transversely in one point. It is well known (and easy to prove)
that under such conditions the given representation is in fact
a (1,2) decomposition of $K(0,1)$ relative to either $S_0$ or
$S_1$. A (1,2) decomposition for the knot $K(p,q)$ can be
obtained by Dehn twisting the knot $K(0,1)$ $p$-times along the
disk $D_0$ in Fig.~\ref{fig11e} (following our previous
convention) and varying the number $q$ of full twists on the
top pair of strands that run between $S_0$ and $S_1$ (with
$q>0$ corresponding to $q$ positive full twists on the
strands).

On the other hand, it is not easy to see that most of the knots
$K(p,q)$ are not $(1,1)$-knots. Our main result is the
following:

\begin{thm}\label{main}
The $(1,2)$-knot $K(p,q)$ is trivial iff $(p,q)=(0,0),(0,1)$,
and a torus knot iff $p=0$ and $q\neq -1,0$ (with
$K(0,q)=T(2,2q-1)$) or $(p,q)=(-1,1)$ (with $K(-1,1)=T(5,8))$;
in all other cases,
\ben
\item[(a)] $K(p,q)$ is a hyperbolic tunnel number one knot
which is not 2-bridge;

\item[(b)]
$K(p,q)$ bounds an unknotted Seifert Klein bottle $P(p,q)$ of
boundary slope $r=4q-36p$;

\item[(c)]
$K(p,q)$ is a $(1,1)$-knot iff $(p,q)$ is a pair of the form
$(p,1)$, $(p,2)$, $(1,q)$, or $(-1,0)$;

\item[(d)]
with the exception of $K(-1,2)(r)=S^2(2,2,3)$,
$K(-2,1)(r)=S^2(2,2,7)$, and $K(p,0)(r)=S^2(2,2, |6p-1|)$, the
manifold $K(p,q)(r)$ is irreducible and toroidal.
\een

In particular, there are infinitely many hyperbolic
$(1,2)$-knots with crosscap number two and tunnel number one
which admit no $(1,1)$ decompositions.
\end{thm}

\bpf
We have already seen that each $K(p,q)$ is a (1,2)-knot; all
other claims follow from Lemmas~\ref{hyper2}, \ref{toro}, and
\ref{one}.
\epf

We remark that the family of knots $K(p,q)$ is one of the
simplest families that can be obtained following our method of
construction, which is quite general. The paper is organized as
follows. In Section~\ref{sec2} we provide definitions and
background results, and develop several specific properties of
circles embedded in the boundary of a genus two handlebody, in
both algebraic and topological versions, which will be needed
in later sections. Section~\ref{sec3} contains, among many
other miscellaneous results, the criteria used to determine if
a crosscap number two knot with an unknotted Seifert Klein
bottle has tunnel number one, and if so whether or not it
admits a (1,1) decomposition; such criteria are given in
decidable algebraic terms involving primitive or power words in
a rank two free group. Finally, in Section~\ref{sec4} we apply
these criteria to prove Lemmas~\ref{hyper2}, \ref{toro} and
\ref{one}, which establish the properties of the family of
knots $K(p,q)$ given in Theorem~\ref{main}.

\section{Preliminaries}\label{sec2}

\subsection{Once-punctured Klein bottles}\label{kbottles}

Let $P$ denote a once-punctured Klein bottle. Any circle
embedded in $P$ is, up to isotopy, of one of the following
types (cf \cite[Section 2]{valdez6}):
\ben
\item[(i)] a {\it meridian circle $m$:} this is an orientation
preserving circle which cuts $P$ into a pair of pants;

\item[(ii)] a {\it center $c$:} this is an orientation reversing
circle whose regular neighborhood in $P$ is a Moebius band;

\item[(iii)] a {\it longitude $\ell$:} this is an orientation preserving
circle which separates $P$ into two components, one of which is a
Moebius band.
\een

The meridian circle of $P$ is unique, while there are
infinitely many isotopy classes of center and longitude circles
(cf \cite[Lemma 3.1]{valdez7}). This contrasts with the
situation in a {\it closed\/} Klein bottle, where up to isotopy
there is only one longitude circle and two center circles.

Denote by $P\widetilde{\times} I$ the orientable twisted
$I$-bundle over $P$, where $I=[0,1]$. $P\widetilde{\times} I$
is a genus two handlebody; the pair $(P\wt{\times}I,\partial
P)$ can be seen in Fig.~\ref{fig13}, up to homeomorphism. In
particular, if $N(P)$ is the regular neighborhood of a
once-punctured Klein bottle which is properly embedded in an
orientable 3--manifold with boundary, then $(N(P),P)\approx
(P\widetilde{\times} I,P\widetilde{\times}
\frac{1}{2})$, where $\approx$ denotes homeomorphism.

\subsection{Lifts of meridian, center, and longitude circles}

Let $T_P$ be the twice punctured torus  $\partial
N(P)\setminus\intr N(\partial P)\subset \partial N(P)$. For any
meridian circle $m$ or longitude circle $\ell$ of $P$, the
restriction of the $I$-bundle $N(P)$ to $m,\ell$ is a fibered
annulus $A(m),A(\ell)$, respectively, properly embedded in
$N(P)$, which intersects $P$ transversely in $m,\ell$,
respectively; if $c$ is a center circle of $P$, the restriction
of $N(P)$ to $c$ is a fibered Moebius band $B(c)$ which
intersects $P$ transversely in $c$. Notice that $A(m),A(\ell),$
and $B(c)$ are all unique up to isotopy in $N(P)$ (ie, any
annulus $(A,\partial A)\subset (N(P),T_P)$ intersecting $P$ in
$m$ is isotopic to $A(m)$, etc.), and that the boundary circles
$\partial A(m),\partial A(\ell),\partial B(c)$ may all be
assumed to lie in $T_P$. We call the circles $\partial
A(m),\partial A(\ell),\partial B(c)$, respectively, the {\it
lifts} of $m,\ell,c$ to $T_P$; thus $m$ has two distinct
nonparallel lifts $m_0\sqcup m_1=\partial A(m)$, while $c$ has
a unique lift. Each longitude $\ell$ also has a unique lift, of
which $\partial A(\ell)$ gives two parallel copies. In fact, if
$\ell$ splits off a Moebius band from $P$ with center
$c_{\ell}$, then the lifts of $c_{\ell}$ and $\ell$ are
isotopic in $T_P$: for $A(\ell)$ is isotopic to the frontier of
the regular neighborhood in $N(P)$ of the Moebius band
$B(c_{\ell})$. Thus, the set of lifts of centers of $P$
coincides with the set of lifts of longitudes of $P$.

\subsection{Seifert Klein bottles}

Let $P\subset X_K$ be a Seifert Klein bottle for a knot
$K\subset S^3$, and let $N(P)\approx P\wt{\times}I$ be a small
regular neighborhood of $P$ in $X_K$. We define the {\it
exterior of $P$ in $S^3$} as the manifold
$X(P)=S^3\setminus\intr N(P)$; we thus have
\begin{equation}\label{basic}
S^3=N(P)\cup_{\partial}X(P)
\end{equation}
with $\partial P\subset\partial N(P)=\partial X(P)$. We will
identify the twice punctured torus $T_P\subset\partial N(P)$
with the frontier in $X_K$ of $N(P)$, so that $T_P\subset N(P)
\cap\partial X(P)$.

Given that $K$ and $\partial P$ are isotopic in $S^3$, the
translation of properties of $K\subset S^3$ or $P\subset X_K$
into properties involving the decomposition given in
(\ref{basic}) can be easily carried out. For instance, it is
easy to see that $P$ is unknotted in $S^3$ iff $X(P)$ is a
handlebody.

\subsection{Companion annuli and multiplicity}\label{companion}

Let $\mc{M}$ be an orientable, irreducible, and geometrically
atoroidal 3--manifold with connected boundary, and let $\gamma$
be a circle embedded in $\partial \mc{M}$ which is nontrivial
(ie, it does not bound a disk) in $\mc{M}$.

Let $A$ be an annular regular neighborhood of $\gamma$ in
$\partial\mc{M}$, and $A'$ a properly embedded separating
annulus in $\mc{M}$ with $\partial A'=\partial A$. We say that
$A'$ is a {\it  companion annulus for $\gamma$ in $\mc{M}$} if
$A'$ is not parallel into $\partial\mc{M}$. It follows from
\cite[Lemma 5.1]{valdez7} that the region cobounded in $\mc{M}$
by $A'$ and the annular neighborhood $A$ of $\gamma$ is a solid
torus, the {\it companion solid torus of $\gamma$ in $M$}, and
that a companion annulus and a companion solid torus for
$\gamma$ are unique up to isotopy. We define the {\it
multiplicity $\mu(\gamma)$ of $\gamma$} in $\mc{M}$ to be $1$
if $\gamma$ has no companion annuli, and as the number of times
$\gamma$ runs around its companion solid torus when $\gamma$
has a companion annulus. Thus $\gamma$ has a companion annulus
in $\mc{M}$ iff $\mu(\gamma)\geq 2$.

Multiplicities of circles in the case where $\mc{M}$ is a genus
two handlebody will be of particular interest in later
developments. So let $H$ be a genus two handlebody; we shall
see that the fact that $\pi_1(H)$ (rel some base point) is a
free group on two generators allows for a simple interpretation
of multiplicities of circles in $\partial H$ in purely
algebraic terms. We will need the following general
definitions.

Let $F_2$ denote the free group on the two generators; if free
generators (ie, a basis) $x,y$ for $F_2$ are given, so that
$F_2=\langle x,y \ | -\rangle$, we may refer to the elements of
$F_2$ as words in $x$ and $y$. For convenience, we will denote
the inverse $u^{-1}$ of $u$ by $\overline{u}$, and by $[u,v]$
the commutator $uv\,\overline{u}\,\overline{v}$ of any two
elements $u,v$ of $F_2$. A word $u\in F_2$ is {\it primitive}
if there is $v\in F_2$ such that $\{u,v\}$ is a basis of $F_2$,
and that $u$ is a {\it power} if there is a nontrivial element
$w\in F_2$ and an integer $n\geq 2$ such that $u=w^n$. We write
$u\equiv v$ for $u,v\in F_2$ if $u=\ovc\,v^{\ve}\,c$ for some
$c\in F_2$ and $\ve\in\{1,-1\}$. A word $u\in\grp{x,y
\ | \ -}$ is said to be {\it cyclically reduced} if, for
$\ve=\pm 1$ and $z\in\{x,y\}$, the pair of symbols
$z^{\ve},z^{-\ve}$ do not occur consecutively in $u$ nor $u$
simultaneously begins with $z^{\ve}$ and ends with $z^{-\ve}$;
notice $u\equiv v$ whenever $v$ is a cyclic reduction of $u$ or
$\overline{u}$.

We will see in Lemma~\ref{comp} below that multiplicities of
circles in $\partial H$ can be characterized in terms of
primitive or power elements in $\pi_1(H)$. A complete
characterization of primitive words in $F_2=\grp{x,y
\ | \ -}$ can be found in \cite{fico1}; for our purposes, the
following  partial characterization of such primitive words
(originally given in \cite{cohen}), which easily extends to
words that are powers of primitive elements, will suffice:

\begin{lem}(\cite{cohen},\cite{fico1})\label{cohen}
If an element $u$ in the free group $F_2=\grp{x,y \ |
\ -}$ is primitive or a power of a primitive then there is a basis
$\{a,b\}\subset\{x,\ovx,y,\ovy\}$ of $F_2$ and an integer
$n\geq 1$ such that either $u\equiv ab^{n-1}$, $u\equiv
a^{n+1}$, or $u\equiv ab^{m_1}\cdots ab^{m_k}$, $k\geq 2$, for
some integers $\{m_1,\dots,m_k\}=\{n,n+1\}$.\hfill\qed
\end{lem}

The next result establishes some useful equivalences in $F_2$.

\begin{lem}\label{free}
Let  $\{u,v\}$ and $\{a,b\}$ be any two bases of $F_2$ and
$m,n\geq 1$ any two integers.
\ben
\ita
The identities $[\overline{w},w']\equiv [w,w']\equiv [w',w]$
hold for any $w,w'\in F_2$.

\itb
If $w\in F_2$, then $\{u,w\}$ is a basis for $F_2$ iff $w=u^k
v^{\ve} u^l$ for some integers $k,l$ and $\ve\in\{-1,1\}$.

\itc
If $u^m\equiv a^m$ then $[u^m,v]\equiv [a^m,b]$.

\itd
If $u^m\equiv a^m$ and $v^n\equiv b^n$ then $[u^m,v^n]\equiv
[a^m,b^n]$.
\een
\end{lem}

\bpf
Part (a) follows by direct computation, while the result in (b)
is well known (cf \cite[Section 3.5, problem 3]{mks}).

For part (c) we have that
$u^m=\ovc\,a^{m\ve}\,c=(\ovc\,a^{\ve}\,c)^m$ for some $c\in
F_2$ and $\ve\in\{-1,1\}$, and hence that $u=\ovc\,a^{\ve}\,c$
(cf \cite[Section 1.4]{mks}). It follows that $cu\ovc=a^{\ve}$
and $cv\ovc$ form a basis for $F_2$, and hence by (b) that
$cv\ovc=a^kb^{\ve'}a^l$ for some integers $k,l$ and
$\ve'\in\{-1,1\}$. Therefore
$$
[u^m,v]=[\ovc a^{m\ve}c,\ovc a^kb^{\ve'}a^lc]\equiv
[a^{m\ve},b^{\ve'}]\equiv [a^m,b].
$$

Similarly, for part (d) we continue to have
$u^m=\ovc\,a^{m\ve}\,c=(\ovc\,a^{\ve}\,c)^m$ and also have
$v^n=\ovd\,b^{n\delta}\,d=(\ovd\,b^{\delta}\,d)^n$ for some
$d\in F_2$ and $\delta\in\{-1,1\}$, whence $u=\ovc\,a^{\ve}\,c$
and $v=\ovd\,b^{\delta}\,d$ hold in $F_2$. It follows that
$cu\ovc=a^{\ve}$ and $cv\ovc=\ove b^{\delta}e$ form a basis of
$F_2$, where $e=d\ovc$. By (b) we must have $\ove
b^{\delta}e=\ova^kb^{\delta}a^k$, and so
$[u^m,v^n]=[\ovc a^{m\ve}c,\ovc\, \ova^kb^{\delta}a^k c]\equiv
[a^{m\ve},b^{\delta}]\equiv [a^m,b^n]$ holds.
\epf

The following sequence of lemmas will establish several
fundamental facts about circles in $\partial H$ which may
represent primitive/power elements in $\pi_1(H)$. For any loop
$\alpha\subset H$, denote by $[\alpha]$ the element of
$\pi_1(H)$ representing $\alpha$ (rel some base point). We will
call any disk properly embedded in $H$ which separates $H$ into
two solid torus components a {\it waist disk} of $H$. Also, for
a 3--manifold $\mc{M}$, we will say that the pair
$(\mc{M},\partial\mc{M})$ is {\it irreducible} if $\mc{M}$ is
irreducible and $\partial\mc{M}$ is incompressible in $\mc{M}$.

\begin{lem}\label{comp}
Let $H$ be a handlebody of genus two and $\gamma,\gamma'$ be
disjoint circles embedded in $\partial H$ which are nontrivial
in $H$.
\ben
\ita
$\partial H\setminus\gamma$ compresses in $H$ iff $[\gamma]$ is
primitive or a power in $\pi_1(H)$, in which case there is a
waist disk of $H$ which is disjoint from $\gamma$. In
particular, if $[\gamma]$ is a power in $\pi_1(H)$ then
$[\gamma]$ is a power of a primitive element.

\itb
$\gamma$ has multiplicity $n\geq 2$ in $H$ iff $[\gamma]$ is
the $n$th power of some primitive element in $\pi_1(H)$;
moreover, if $[\gamma]=\lambda^n$ for some integer $n\geq 1$
and some primitive element $\lambda\in\pi_1(H)$ then
$n=\mu([\gamma])$.

\itc
Suppose $\gamma$ and $\gamma'$ are not parallel in $\partial
H$. If $[\gamma]$ is primitive or a power in $\pi_1(H)$ and
$[\gamma']$ is conjugate to $[\gamma]$ in $\pi_1(H)$, then
$\gamma$ and $\gamma'$ cobound a nonseparating annulus in $H$.

\itd
Suppose $\gamma'$ does not cobound an annulus with $\gamma$ in
$H$ and that, in $\pi_1(H)$, either $[\gamma]$ is primitive or
a power while $[\gamma']$ is a power. Then there is a waist
disk $D$ of $H$ which separates $\gamma$ and $\gamma'$.
\een
\end{lem}

\bpf
Parts (a) and (b) follow from the argument used in the proof of
\cite[Theorem 4.1]{cassongor}, which deals with roots in the
fundamental group of a compression body; we prove them here in
the context of handlebodies of genus two for the convenience of
the reader.

Suppose $\partial H\setminus\gamma$ compresses in $H$ along a
disk $D$. If $D$ is nonseparating then there is an embedded
circle $\alpha\subset\partial H\setminus\gamma$ which
intersects $\partial D$ transversely in a single point, and
hence the frontier of a small regular neighborhood in $H$ of
$D\cup\alpha$ is a waist disk in $H$ which compresses $\partial
H\setminus\gamma$; we may thus assume $D$ is a waist disk of
$H$. Cutting $H$ along $D$ produces two solid tori components,
one of which, say $V$, contains $\gamma$ in its boundary; given
that $\gamma$ is nontrivial in $H$, it follows that $\gamma$
must be a nonseparating circle in $\partial H$. By Van Kampen's
theorem, if $\beta$ is a core of $V$ then $[\beta]$ is
primitive in $\pi_1(H)$ and $[\gamma]=[\beta]^k$ for some
integer $k\neq 0$. Thus $[\gamma]$ is either primitive or a
power of a primitive in $\pi_1(H)$; moreover, it is not hard to
see that if $|k|\geq 2$ then $V$ is the companion solid torus
of $\gamma$ in $H$, so $\mu(\gamma)=|k|$. This proves one
direction of part (a).

Conversely, let $M$ be the 3--manifold obtained by adding a
2-handle to $H$ along $\gamma$. If, in $\pi_1(H)$, $[\gamma]$
is primitive then $\pi_1(M)=\mZ$ and so $M$ is a solid torus,
while if $[\gamma]$ is a power then $\pi_1(M)$ has nontrivial
torsion by \cite[Theorems N3 and 4.12]{mks} and hence $M$ is
reducible (cf \cite[Theorem 9.8]{hempel}). Therefore the pair
$(M,\partial M)$ is not irreducible, so the surface $\partial
H\setminus\gamma$ compresses in $H$ by the 2-handle addition
theorem (cf \cite{cassongor}), and hence by the above argument
$[\gamma]$ is a primitive or a power of a primitive in
$\pi_1(H)$. Thus (a) holds.

For part (b), assume $\gamma$ has multiplicity $n\geq 2$; that
is, for $A\subset\partial H$ an annular neighborhood of
$\gamma$ and $A'$ a companion annulus of $\gamma$ with
$\partial A'=\partial A$, $A$ and $A'$ cobound a solid torus
$V\subset H$ such that $\gamma$ runs $n$ times around $V$.
Since $\gamma$ is nontrivial in $H$, and $A'$ separates $H$, it
follows that $A'$ boundary compresses in $H$ into a nontrivial
separating compression disk of $\partial H\setminus\gamma$;
hence, by the first part of the argument for (a), $\gamma$ is
the $n$th power of a primitive element of $\pi_1(H)$.
Conversely, suppose $[\gamma]=\lambda^n$ for some integer
$n\geq 1$ and some primitive element $\lambda\in\pi_1(H)$. By
the first part of the argument for (a), there is a loop $\beta$
in $H$ with $[\beta]$ primitive in $\pi_1(H)$ and
$[\gamma]=[\beta]^{\mu(\gamma)}$. Therefore the abelianization
of $\pi_1(H)/
\langle[\gamma]\rangle$ is isomorphic to both $\mathbb{Z}
\oplus\mathbb{Z}_n$ and $\mathbb{Z}
\oplus\mathbb{Z}_{\mu(\gamma)}$, so $n=\mu(\gamma)$. Thus (b)
follows.

In parts (c) and (d), let $F=\partial H\setminus\gamma$, so
that $\gamma'\subset F$; observe that $F$ compresses in $H$ by
(a), given that $[\gamma]$ is either primitive or a power in
$\pi_1(H)$. Let $M$ be the manifold obtained by attaching a
2-handle to $H$ along $\gamma'$, so that $\gamma\subset\partial
M$. In part (d), $[\gamma']$ is a power in $\pi_1(H)$ and so
$M$ is reducible by the argument used in part (a); in part (c),
since $[\gamma]$ is conjugate to $[\gamma']$ in $\pi_1(H)$ but
$\gamma$ and $\gamma'$ are not parallel in $\partial H$, and
$\gamma\subset\partial M$, $[\gamma]$ must be trivial in
$\pi_1(M)$ but nontrivial in $\partial M$, and hence $\gamma$
bounds a nonseparating disk in $M$. Either way the pair
$(M,\partial M\setminus\gamma)$ is not irreducible, so by the
2-handle addition theorem the surface
$F\setminus\gamma'=\partial H\setminus(\gamma\cup\gamma')$
compresses in $H$ along some disk $D\subset H$. In part (c) the
disk $D$ must be nonseparating, so $\gamma,\gamma'$ lie in the
boundary of the solid torus $H\setminus\intr N(D)$ and hence
cobound a nonseparating annulus in $H$; similarly, in part (d)
the disk $D$ must be a waist disk of $H$ which separates
$\gamma$ and $\gamma'$.
\epf

For a Seifert Klein bottle $P$ in a knot exterior $X_K$, recall
that $T_P$ is the twice punctured torus obtained from the
frontier of $N(P)$ in the knot exterior $X_K$, so
$T_P\subset\partial X(P)\setminus\partial P$ and $\intr T_P$,
$X(P)\setminus\partial P$ are homeomorphic surfaces. In this
context, Lemma~\ref{comp}(a) has the following immediate
consequence.

\begin{cor}\label{pi1}
An unknotted Seifert Klein bottle $P$ for a knot $K$ is
$\pi_1$-injective in $X_K$ iff $[\partial P]$ is neither
primitive nor a power of a primitive in $\pi_1(X(P))$.
Specifically, $T_P$ is boundary compressible in $X_K$ iff
$[\partial P]$ is primitive in $\pi_1(X(P))$.
\end{cor}

\bpf
Since $N(P)$ is an $I$-bundle over $P$, we have that $P$ is
$\pi_1$-injective in $N(P)$ and $T_P$ is incompressible in
$N(P)$; hence, by the Dehn's lemma-loop theorem \cite{hempel},
$P$ is $\pi_1$-injective in $X_K$ iff $T_P$ is geometrically
incompressible in $X(P)$. Thus the first part follows from
Lemma~\ref{comp}(a).

Now, given the relationship between the surfaces $T_P$ and
$X(P)\setminus\partial P$, it is not hard to see that the
boundary compressibility of $T_P$ in $X(P)$ is equivalent to
the existence of a properly embedded disk $D$ in $X(P)$ which
intersects the circle  $\partial P\subset\partial X(P)$
transversely in one point; as the latter condition is
equivalent to $[\partial P]$ being primitive in $\pi_1(X(P))$,
the second part of the claim follows.
\epf

The following result gives a simple algebraic way of
determining if the manifold $T\times I$, $T$ a torus, is
obtained by attaching a 2-handle to a genus two handlebody $H$;
though the result is well known, we sketch its proof as
preparation for the argument used in its generalization given
in Lemma~\ref{comp2}, which deals with the case of attaching a
2-handle to a genus two sub-handlebody of $H$.

\begin{lem}\label{ti}
Let $H$ be a genus two handlebody and $\mc{T}$ a closed torus.
Let $\gamma$ be a circle embedded in $\partial H$ and
$M=H\cup_{\gamma}N(D)$ be the manifold obtained by a attaching
a 2-handle $N(D)$ to $H$ along $\gamma$. Then
$M\approx\mc{T}\times I$ iff $[\gamma]\equiv [a,b]$ for some
(and hence any) basis $\{a,b\}$ of $\pi_1(H)$.
\end{lem}

\bpf
Suppose that $\pi_1(M)/ \grp{[\gamma]}=\mZ\oplus\mZ$, so that
$\gamma$ is nontrivial in $H$. If $\partial H\setminus\gamma$
compresses in $H$ then by Lemma~\ref{comp}(a) there is a waist
disk in $H$ disjoint from $\gamma$ and so $M$ is a manifold of
the form $S^1\times D^2\# L$ for $L=S^3$ or a lens space; but
then $\pi_1(M)/ \grp{[\gamma]}\neq\mZ\oplus\mZ$, contradicting
our hypothesis. Thus $\partial H\setminus\gamma$ is
incompressible in $H$ and so the pair $(M,\partial M)$ is
irreducible by the 2-handle addition theorem, hence by
\cite[Theorem 12.10]{hempel} the condition
$M\approx\mc{T}\times I$ is equivalent to the condition
$\pi_1(H)/ \grp{[\gamma]}=\mZ\oplus\mZ$. The lemma follows now
from the fact (due to Nielsen, cf \cite[Section 4.4]{mks})
that, for any word $w$ in $\grp{x,y \ | \ -}$,  $\grp{x,y \ | \
-}/
\grp{w}=\mZ\oplus\mZ$ iff $w\equiv[x,y]$.
\epf

\begin{lem}\label{comp2}
Let $H$ be a handlebody of genus two and
$\gamma_0,\gamma_1,\gamma_2$ be disjoint circles embedded in
$\partial H$ which are nontrivial in $H$; let $\mc{T}$ denote a
closed torus.

\ben
\ita If $A_0\subset H$ is a companion annulus for $\gamma_0$
with core $\alpha_0$ and corresponding companion solid torus
$V_0\subset H$, then $H'=\clo(H\setminus V_0)$ is a genus two
handlebody, and there is a common waist disk $D$ of $H$ and
$H'$ such that
\ben
\item[(i)]
$H'=W_0\cup_{D} W_1$ for some solid tori $W_0,W_1$ in $H'$ with
$D=\partial W_0\cap\partial W_1=W_0\cap W_1$ and
$A_0\subset\partial W_0\setminus D$,

\item[(ii)]
if $\beta_1$ is a core of $W_1$ then
$\{w_0=[\alpha_0],w_1=[\beta_1]\}$ is a basis for $\pi_1(H')$,

\item[(iii)] if $\beta_0$ is a core of $V_0$ then
$\{u=[\beta_0],v=[\beta_1]\}$ is a basis for $\pi_1(H)$, and
the inclusion map $i\co H'\subset H$ induces an injection
$\pi_1(H')\overset{i_*}{\to}\pi_1(H)$ given by $w_0\mapsto
u^{\mu(\gamma_0)}\equiv[\gamma_0]$ and $w_1\mapsto v$,

\item[(iv)] if $H'\cup N(D(\gamma_2))$ is the manifold obtained by
attaching a 2-handle $N(D(\gamma_2))$ along $\gamma_2$, then
$H'\cup N(D(\gamma_2))\approx\mc{T}\times I$ iff
$[\gamma_2]\equiv [a^{\mu(\gamma_0)},b]$ for some (and hence
any) basis $\{a,b\}$ of $\pi_1(H)$ such that $[\gamma_0]\equiv
a^{\mu(\gamma_0)}$.
\een

\itb
Suppose $A_0,A_1\subset H$ are disjoint companion annuli for
$\gamma_0,\gamma_1$, respectively, with corresponding cores
$\alpha_0,\alpha_1$ and companion solid tori $V_0,V_1\subset
H$. Then $H'=\clo(H\setminus (V_0\sqcup V_1))$ is a genus two
handlebody, and there is a common waist disk $D$ of $H$ and
$H'$ such that

\ben
\item[(i)]
$H'=W_0\cup_{D} W_1$ for some solid tori $W_0,W_1$ in $H'$ with
$D=\partial W_0\cap\partial W_1=W_0\cap W_1$,
$A_0\subset\partial W_0\setminus D$, and $A_1\subset\partial
W_1\setminus D$,

\item[(ii)]
$\{w_0=[\alpha_0],w_1=[\alpha_1]\}$ is a basis for $\pi_1(H')$,

\item[(iii)] if $\beta_0,\beta_1$ are cores of $V_0,V_1$, respectively,
then $\{u=[\beta_0],v=[\beta_1]\}$ is a basis for $\pi_1(H)$
and the inclusion map $i\co H'\subset H$ induces an injection
$\pi_1(H')\overset{i_*}{\to}\pi_1(H)$ given by $w_0\mapsto
u^{\mu(\gamma_0)}\equiv [\gamma_0]$ and $w_1\mapsto
v^{\mu(\gamma_1)}\equiv [\gamma_1]$.

\item[(iv)]
if $H'\cup N(D(\gamma_2))$ is the manifold obtained by
attaching a 2-handle $N(D(\gamma_2))$ along $\gamma_2$, then
$H'\cup N(D(\gamma_2))\approx\mc{T}\times I$ iff
$[\gamma_2]\equiv [a^{\mu(\gamma_0)},b^{\mu(\gamma_1)}]$ for
some (and hence any) basis $\{a,b\}$ of $\pi_1(H)$ such that
$[\gamma_0]\equiv a^{\mu(\gamma_0)}$ and $[\gamma_1]\equiv
b^{\mu(\gamma_1)}$.
\een
\een
\end{lem}

\bpf
For part (a), by Lemma~\ref{comp}(a,b), there is a waist disk
$D$ for $H$ such that $H=U_0\cup_D W_1$ for some solid tori
$U_0,W_1$ with $\gamma_0\subset\partial U_0\setminus D$. After
a slight isotopy we may also assume that $A_0\subset U_0$, and
then we may write $U_0=W_0\cup_{A_0} V_0$ for some solid torus
$W_0\subset U_0$. Thus $H'=\clo(H\setminus V_0)=W_0\cup_D W_1$
is a genus two handlebody and (i) holds.

As the circles $\gamma,\alpha\subset\partial V_0$ run
$\mu(\gamma_0)\geq 2$ around $V_0$, and $U_0=W_0\cup_{A_0} V_0$
is a solid torus, it follows that $\alpha_0\subset\partial W_0$
must run once around $W_0$ and hence that $\alpha_0$ is
isotopic to a core of $W_0$; therefore, that (ii) and (iii)
hold follows by Van Kampen's theorem and the fact that
$\alpha_0$ and $\gamma_0$ are isotopic in $V_0$.

For part (a)(iv), we assume as we may that $\partial A_0$ and
$\gamma_2$ are disjoint in $\partial H$, whence
$\gamma_2\subset\partial H'$; we write $[\gamma_2]',[\gamma_2]$
for the elements in $\pi_1(H'),\pi_1(H)$ represented by
$\gamma_2$, respectively, so that
$i_*([\gamma_2]')=[\gamma_2]$. Recall by Lemma~\ref{ti} that
$H'\cup N(D(\gamma_2))\approx\mc{T}\times I$ iff
$[\gamma_2]'\equiv[x,y]$ for some and in fact any basis
$\{x,y\}$ of $\pi_1(H')$. Thus, if $H'\cup
N(D(\gamma_2))\approx\mc{T}\times I$ then
$[\gamma_2]'\equiv[w_0,w_1]$ in $\pi_1(H')$ and hence
$[\gamma_2]=i_*([\gamma_2]')\equiv [u^{\mu(\gamma_0)},v]$ in
$\pi_1(H)$; that $[\gamma_2]\equiv [a^{\mu(\gamma_0)},b]$ holds
for any basis $\{a,b\}$ of $\pi_1(H)$ with $[\gamma_0]\equiv
a^{\mu(\gamma_0)}$ now follows from Lemma~\ref{free}(c).

Suppose now $\{a,b\}$ is any basis of $\pi_1(H)$ such that
$[\gamma_0]\equiv a^{\mu(\gamma_0)}$ and $[\gamma_2]\equiv
[a^{\mu(\gamma_0)},b]$ hold in $\pi_1(H)$; by
Lemma~\ref{free}(c), we then also have that $[\gamma_2]\equiv
[u^{\mu(\gamma_0)},v]$. Observe that $[u^{\mu(\gamma_0)},v]$ is
a cyclically reduced word in $\pi_1(H)=\grp{u,v \ | \ -}$; for
definiteness, we will assume, that $[u^{\mu(\gamma_0)},v]$ is a
cyclic reduction of $[\gamma_2]$ in $\pi_1(H)=\grp{u,v \ |
\ -}$.

Let $W(w_0,w_1)$ be a cyclic reduction of $[\gamma_2]'$ in
$\pi_1(H')=\grp{w_0,w_1 \ | \ -}$. Then, in $\pi_1(H)=\grp{u,v
\ | \ -}$,
$i_*(W(w_0,w_1))=W(i_*(w_0),i_*(w_1))=W\left(u^{\mu(\gamma_0)},v\right)$
is also a cyclically reduced word, which must then be a cyclic
reduction of $[\gamma_2]=i_*([\gamma_2]')$. Thus the words
$i_*(W(w_0,w_1))$ and $[u^{\mu(\gamma_0)},v]$ are identical
except for the cyclic order of their factors. Given that
$i_*(w_0)=u^{\mu(\gamma_0)}$ and $i_*(w_1)=v$, it is not hard
to see that changing the cyclic order of the $w_0,w_1$ factors
in $W(w_0,w_1)$ will produce the identity $i_*(W(w_0,w_1))=
[u^{\mu(\gamma_0)},v]$, so we may assume that such identity
holds. As $i_*([w_0,w_1])=[u^{\mu(\gamma_0)},v]$ holds too, so
that $i_*([w_0,w_1])=i_*(W(w_0,w_1))$, we must have $[w_0,w_1]=
W(w_0,w_1)\equiv[\gamma_2]'$ in $\pi_1(H')$. Thus (a)(iv)
holds.

The proof for part (b) is similar: by Lemma~\ref{comp}(d),
there is a waist disk $D$ for $H$ such that $H=U_0\cup_D U_1$
for some solid tori $U_0,U_1$ with $\gamma_0\subset\partial
U_0\setminus D$ and $\gamma_1\subset\partial U_1\setminus D$,
and we may also assume that $A_0\subset U_0, A_1\subset U_1$.
Thus $U_0=W_0\cup_{A_0} V_0$ and $U_1=W_1\cup_{A_1} V_1$ for
some solid torus $W_0\subset U_0, W_1\subset U_1$, so
$H'=\clo(H\setminus (V_0\cup V_1))=W_0\cup_D W_1$ is a genus
two handlebody and (i) holds. As before, $\alpha_0,\alpha_1$
are isotopic to cores of $W_0,W_1$, respectively, so that (ii)
and (iii) hold follows by Van Kampen's theorem. The proof of
(b)(iv) follows the same argument as that of (a)(iv), using
Lemma~\ref{free}(d) instead of Lemma~\ref{free}(c) to deduce
$[\gamma_2]\equiv [u^{\mu(\gamma_0)},v^{\mu(\gamma_1)}]$ from
$[\gamma_2]\equiv [a^{\mu(\gamma_0)},b^{\mu(\gamma_1)}]$.
\epf

\section{Crosscap number two knots}\label{sec3}

In this section we establish necessary and sufficient
conditions for a crosscap number two hyperbolic knot to admit a
(1,1) decomposition. We also establish miscellaneous results
that can be used to detect when a crosscap number two knot with
an unknotted Seifert Klein bottle has tunnel number one or is
hyperbolic, as well as means of identifying and constructing
the lifts of meridians, centers, and longitudes of a Seifert
Klein bottle.

\subsection{(1,1) Decompositions}

The following construction of a special family of Seifert Klein
bottles in $S^3$, denoted by $P(c_0^{\varepsilon_0},
c_1^{\varepsilon_1},R)$, with boundary a $(1,1)$-knot, is taken
from \cite[Section 1]{valdez8}. Let $T$ be an unknotted (ie,
Heegaard) torus embedded in $S^3$; we identify a small regular
neighborhood of $T$ in $S^3$ with a product $T\times I$, where
$I=[0,1]$. Thus, there are unknotted solid tori $V_0,V_1\subset
S^3$ such that
\begin{equation}\label{eqn0}
S^3=V_0\,\cup_{\partial V_0=T\times\{0\}}\, T\times I\,
\cup_{\partial V_1=T\times\{1\}}\, V_1.
\end{equation}
We say that an arc $\gamma$ embedded in $T\times I$ is {\it
monotone} if the natural projection map $T\times I\rightarrow
I$ is monotone on $\gamma$. We may further assume that $T\times
I$ lies within a slightly larger embedding of the form
$T\times[-\delta,1+\delta]$, for some small $\delta>0$.

For $i=0,1$, let $c_i$ be a circle nontrivially embedded in
$T\times\{i\}$. Let $R$ be a rectangle properly embedded in
$T\times I$ with one boundary side along $c_0$ and the opposite
side along $c_1$, such that $R\cap (T\times [0,\delta]\cup
T\times [1-\delta,1])\subset c_0\times [0,\delta]\cup c_1\times
[1-\delta,1]$ and some core $\beta\subset T\times I$ of $R$ is
monotone. The union of $R$ with the annuli
$\mc{A}_i=c_i\times[i-\delta,i+\delta]$, $i=0,1$ is then a pair
of pants; giving one half-twist relative to $T\times\{i\}$ to
each annulus piece $\mc{A}_i$, away from $R$, produces a once
punctured Klein bottle $P(c_0^{\varepsilon_0},
c_1^{\varepsilon_1},R)$, where $\varepsilon_i\in\{+,-\}$ and
the notation $c_i^{\varepsilon_i}$ stands for one of the two
possible half-twists that can be performed on the annulus
$\mc{A}_i$ (see Fig.~\ref{fig07}). We remark that in
\cite{valdez8} the knot $\partial P(c_0^{\varepsilon_0},
c_1^{\varepsilon_1},R)$ is denoted by $K(c_0^{*}, c_1^{*},R)$.

\begin{figure}
\psfrag{qI}{$\{q\}\times I$}
\psfrag{Ti}{$T\times\{i\}$}
\psfrag{half-twist}{half-twist on $\mc{A}_i$}
\psfrag{A_2}{$\mc{A}_i$}
\psfrag{A_0}{$A_i$}
\psfrag{A_1}{$A_i'$}
\psfrag{V_0}{$V_i$}
\psfrag{R}{$R$}
\psfrag{TI}{$T\times I$}
\psfrag{PP}{$\partial P(c_0^{\varepsilon_0},
c_1^{\varepsilon_1},R)$}
\psfrag{ci}{$c_i$}
\psfrag{cci}{$\wh{c}_i$}
\Figw{6in}{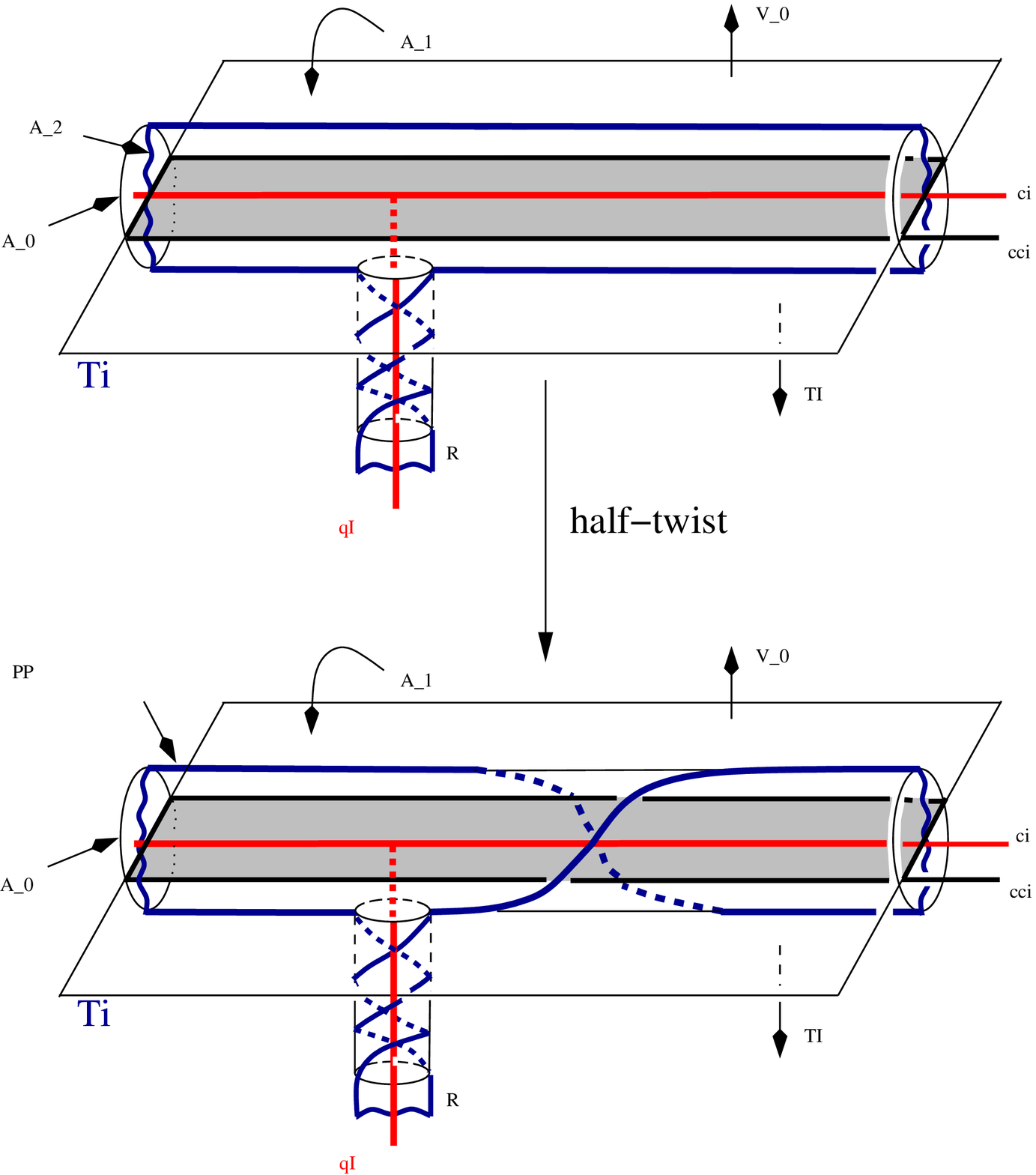}{Construction of the Seifert
Klein bottle $P(c_0^{\varepsilon_0},
c_1^{\varepsilon_1},R)$. }{fig07}
\end{figure}

We say that $P(c_0^{\varepsilon_0}, c_1^{\varepsilon_1},R)$ is
in {\it vertical position} if some monotone core $\beta\subset
R$ is a fiber $\{q\}\times I$ of $T\times I$; in such case,
$P(c_0^{\varepsilon_0}, c_1^{\varepsilon_1},R)$ can be isotoped
so as to properly embed in some regular neighborhood of
$c_0\cup (\{q\}\times I)\cup c_1$ in $S^3$. The next result
states that vertical position is always attainable for any
surface of the form $P(c_0^{\varepsilon_0},
c_1^{\varepsilon_1},R)$.

\begin{lem}\label{vertical}
Any once-punctured Klein bottle of the form
$P(c_0^{\varepsilon_0}, c_1^{\varepsilon_1},R)$ can be isotoped
into vertical position.
\end{lem}

\bpf
Consider an arbitrary once-punctured Klein bottle of the form
$P(c_0^{\varepsilon_0}, c_1^{\varepsilon_1},R)$ with $\beta$ a
monotone core of $R$. We claim that the arc $\beta$ is isotopic
in $T\times I$ to some (and hence any) fiber $\{q\}\times I$,
$q\in T$; in such case, the isotopy that moves $\beta$ onto
some fiber $\{p\}\times I$ can be extended to an isotopy of
$c_0\cup\beta\cup c_1$ in $T\times I$, and then further
extended to an isotopy of $S^3$ that puts
$P(c_0^{\varepsilon_0}, c_1^{\varepsilon_1},R)$ in vertical
position.

Given $0\leq x\leq y\leq 1$ and any monotone arc $\gamma\subset
T\times I$, denote the point $\gamma\cap (T\times\{x\})$ by
$\gamma_x$ and the arc $\gamma\cap (T\times[x,y])$ by
$\gamma_{[x,y]}$. Since $\beta$ is monotone in $T\times I$,
there is a sufficiently large integer $n>0$ such that, for each
integer $1\leq k\leq n$, the arc $\beta_{[(k-1)/n,k/n]}$ is
isotopic, rel $\beta_{k/n}$ (ie, fixing the set
$\{\beta_{k/n}\}$), to the arc $\{\beta_{k/n}\}\times
[(k-1)/n,k/n]$  in $T\times [(k-1)/n,k/n]$.

Isotope the arc $\beta_{[0,1/n]}$ onto the arc
$\{\beta_{1/n}\}\times [0,1/n]$ rel $\beta_{1/n}$, and let
$\beta^{(1)}$ be the union of the arcs $\{\beta_{1/n}\}\times
[0,1/n]$ and $\beta_{[1/n,1]}$; clearly, $\beta^{(1)}$ and $\beta$
are isotopic in $T\times I$ rel $\beta_{[1/n,1]}$. Now isotope the
arc $\beta^{(1)}_{[1/n,2/n]} \ (=\beta_{[1/n,2/n]})$ onto the arc
$\{\beta_{2/n}\}\times [1/n,2/n]$ in $T\times[1/n,2/n]$ rel
$\beta_{2/n}$; this isotopy easily extends to an isotopy of the arcs
$\beta^{(1)}_{[0,2/n]}$  and $\{\beta^{(1)}_{2/n}\}\times [0,2/n]$
in $T\times [0,2/n]$ rel $\beta^{(1)}_{2/n}$, and produces the arc
$\beta^{(2)}=\{\beta^{(1)}_{2/n}\}\times [0,2/n]\cup
\beta^{(1)}_{[2/n,1]}$ in $T\times I$, isotopic to $\beta^{(1)}$
rel $\beta^{(1)}_{[2/n,1]}=\beta_{[2/n,1]}$. Continuing the process
in this fashion, the claim follows by induction, with $\beta$
isotopic to $\beta^{(n)}=\{\beta_{1}\}\times [0,1]$ in $T\times I$
rel $\beta_{1}$.
\epf

\begin{assumpt}\label{as}
In light of Lemma~\ref{vertical}, any Seifert Klein bottle of
the form $\mc{P}=P(c_0^{\varepsilon_0}, c_1^{\varepsilon_1},R)$
constructed relative to an unknotted $T\times I\subset S^3$
will be assumed to be in vertical position relative to $T\times
I$. In particular, we may always assume that
$N(\mc{P})=N(c_0\cup\,\{q\}\times I\, \cup c_1)$ for some point
$q\in T$ and that $T\times I\setminus \intr\, N(\mc{P})$ is
isotopic to $T_0\times I$ for the once punctured torus
$T_0=T\setminus\intr N(q)\subset T$ (see Fig.~\ref{fig07}).
\end{assumpt}

Now let $K$ be a knot in $S^3$ spanning a once-punctured Klein
bottle $\mc{P}$. If $c_0,c_1$ are two disjoint center circles
of $\mc{P}$, we say that {\it $K$ admits a
$\{\mc{P},c_0,c_1\}$-structure} if $\mc{P}$ is isotopic to some
once-punctured Klein bottle of the form $P(c_0^{\varepsilon_0},
c_1^{\varepsilon_1},R)$. In this context, \cite[Theorem
1.1]{valdez8} can be restated, in the case of hyperbolic non
2-bridge knots, as follows:

\begin{thm}{\cite{valdez8}}\label{cafest}
Let $K$ be a hyperbolic knot in $S^3$ which is not 2-bridge and
bounds a Seifert Klein bottle $\mc{P}$. Then $K$ has a $(1,1)$
decomposition iff $K$ admits a $\{\mc{P},c_0,c_1\}$-structure
for some pair of disjoint centers $c_0,c_1$ of $\mc{P}$. \hfill
\qed
\end{thm}

\subsection{Twisted lifts of centers}\label{twisted}
Let $\mc{P}\subset X_K$ be any Seifert Klein bottle for a knot
$K\subset S^3$, and let $N(\mc{P})\subset X_K$ be its regular
neighborhood. If $c_0,c_1$ are disjoint centers of $\mc{P}$,
then, up to isotopy, there is a unique arc $\alpha$ properly
embedded in $\mc{P}$ which separates $c_0$ from $c_1$, and
which gives rise (via the $I$-bundle structure of $N(\mc{P})$)
to a waist disk $D\subset N(\mc{P})$ with $D\cap\mc{P}=\alpha$,
which cuts $N(\mc{P})$ into two solid tori $W_0,W_1$ with
$c_i\subset W_i$ (see Fig.~\ref{fig13}). Notice that such a
waist disk is also unique up to isotopy, and that
$B_i=\mc{P}\cap W_i$ is a Moebius band for $i=1,2$. Performing
one half-twist to $B_i$ in $W_i$ produces an annulus in $W_i$;
there are two ways of half-twisting $B_i$, and each way
produces an annulus, say with boundary slope $\wh{c}_i$ or
$\wh{c}\,'_i\subset\partial W_i\setminus D$, respectively. Each
of the circles $\wh{c}_i,\wh{c}\,'_i$ runs once around $W_i$
and intersects $\partial\mc{P}$ transversely in one point. We
call the circles $\wh{c}_i,\wh{c}\,'_i$ the {\it twisted lifts}
of the center $c_i\subset\mc{P}$ to $\partial N(\mc{P})$.

\begin{figure}
\psfrag{NP}{$N(\mc{P})$}
\psfrag{E}{$\partial \mc{P}$}
\psfrag{P}{$\mc{P}$}
\psfrag{c0}{$c_0$}
\psfrag{c1}{$c_1$}
\psfrag{alpha}{$\alpha$}
\psfrag{cc0}{$\wh{c}_0$}
\psfrag{cc1}{$\wh{c}_1$}
\psfrag{V0}{$W_0$}
\psfrag{V1}{$W_1$}
\psfrag{D}{$D$}
\Figw{6in}{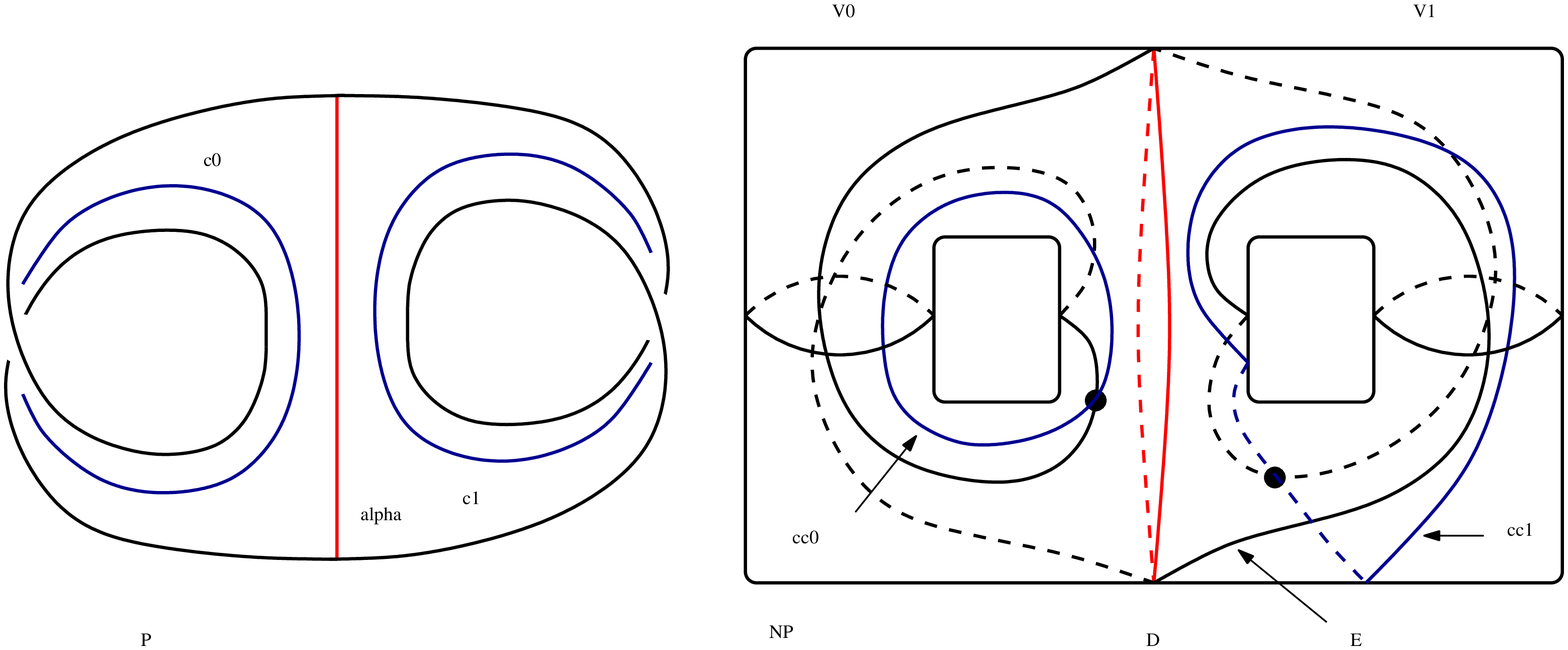}{A pair of disjoint centers $c_0,c_1$ in $\mc{P}$
and some twisted lifts $\wh{c}_0,\wh{c}_1\subset T_{\mc{P}}$.}{fig13}
\end{figure}

In the particular case where $\mc{P}=P(c_0^{\varepsilon_0},
c_1^{\varepsilon_1},R)$ for some disjoint centers
$c_0,c_1\subset\mc{P}$, relative to some unknotted $T\times
I\subset S^3$, so that $N(\mc{P})=N(c_0\cup (\{q\}\times I)\cup
c_1)$, it follows that any circle of intersection $\wh{c}_i$
between $\partial N(c_0\cup (\{q\}\times I)\cup c_1)$ and
$T\times \{i\}$ is a twisted lift of $c_i$ (see
Fig.~\ref{fig07}), which we call the {\it induced twisted lift
of $c_i$}.

Back in the general case, assume further that $\mc{P}$ has
atoroidal exterior $X(\mc{P})\subset S^3$; as $X(\mc{P})$ is
irreducible, companion solid tori and multiplicities of circles
on $\partial X(\mc{P})\subset X(\mc{P})$ are thus defined. So,
for any pair $c_0,c_1$ of disjoint centers of $\mc{P}$ with
twisted lifts $\wh{c}_0,\wh{c}_1\subset
T_{\mc{P}}\subset\partial X(\mc{P})$, define
$V(\wh{c}_i)\subset X(\mc{P})$ as the companion solid torus of
$\wh{c}_i$ if $\mu(\wh{c}_i)\geq 2$ in $X(P)$, with
$A(\wh{c}_i)\subset\partial W_i$ the annular neighborhood of
$\wh{c}_i\subset\partial W_i$ such that $V(\wh{c}_i)\cap
N(\mc{P})= A(\wh{c}_i)$, and otherwise set
$V(\wh{c}_i)=\emptyset=A(\wh{c}_i)$; we now construct the
manifold
\begin{equation}
M(\mc{P},\wh{c}_0,\wh{c}_1)=\clo\left(X(\mc{P})\setminus
\left(V(\wh{c}_0)\cup V(\wh{c}_0)\right)\right)\subset S^3 .
\end{equation}
Let $D$ be a waist disk of $N(\mc{P})$ that separates
$\wh{c}_0$ and $\wh{c}_1$ (see Fig.~\ref{fig13}); $D$ is unique
up to isotopy, and can always be chosen such that
\begin{equation}\label{D}
\partial D\subset\partial M(\mc{P},\wh{c}_0,\wh{c}_1).
\end{equation}
Let $\clo(N(\mc{P})\setminus N(D))=W_0\sqcup W_1$, where
$W_0,W_1$ are solid tori with $\wh{c}_0\subset \partial W_0$
and $\wh{c}_1\subset\partial W_1$. Then
$U_0=W_0\cup_{A(\wh{c}_0)} V(\wh{c}_0)$ and
$U_1=W_1\cup_{A(\wh{c}_1)} V(\wh{c}_1)$ are solid tori since,
whenever present, each annulus $A(\wh{c}_i)$ runs once around
$W_i$ (see Fig.~\ref{fig13}); since, in light of (\ref{basic}),
we have
\begin{equation}\label{eqn3}
S^3= U_0\cup_{\partial} (M(\mc{P},\wh{c}_0,\wh{c}_1)\cup N(D))
\cup_{\partial} U_1,
\end{equation}
it follows that $M(\mc{P},\wh{c}_0,\wh{c}_1)\cup N(D)$ is the
exterior in $S^3$ of the link formed by cores of the solid tori
$U_0$ and $U_1$.

The manifolds $M(\mc{P},\wh{c}_0,\wh{c}_1)$ and
$M(\mc{P},\wh{c}_0,\wh{c}_1)\cup N(D)$ can be readily
identified whenever $\mc{P}$ is of the form
$P(c_0^{\varepsilon_0}, c_1^{\varepsilon_1},R)$.

\begin{lem}\label{mp}
Suppose $\mc{P}=P(c_0^{\varepsilon_0}, c_1^{\varepsilon_1},R)$
for some disjoint centers $c_0,c_1\subset\mc{P}$, relative to
some unknotted $T\times I\subset S^3$; let $\wh{c}_0,\wh{c}_1$
be the induced twisted lifts, respectively, and let $D\subset
N(\mc{P})$ be the unique waist disk that separates $\wh{c}_0$
and $\wh{c}_1$. If $X(\mc{P})$ is atoroidal then there is a
once punctured torus $T_0\subset T$ such that, in $S^3$,
$M(\mc{P},\wh{c}_0,\wh{c}_1)$ is isotopic to $T_0\times I$ and
$M(\mc{P},\wh{c}_0,\wh{c}_1)\cup N(D)$ is isotopic to $T\times
I$, with the isotopy carrying the circle $\partial
D\subset\partial M(\mc{P},\wh{c}_0,\wh{c}_1)$ to a circle in
$\partial (T_0\times I)$ isotopic to $(\partial
T_0)\times\{0\}$.
\end{lem}

\bpf
Recall from (\ref{eqn0}) that $S^3=V_0\cup T\times I\cup V_1$;
also, by our Assumption~\ref{as}, we may assume that
$N(\mc{P})=N(c_0\cup\,\{q\}\times I\,\cup c_1)$ for some $q\in
T$, so that $T\times I\setminus\intr\, N(\mc{P})$ is isotopic
to $T_0\times I$ for the once punctured torus
$T_0=T\setminus\intr N(q)\subset T$.

For $i=0,1$, consider the annuli $A_i\subset\partial X(\mc{P})$
and $A'_i\subset T\times\{i\}$ indicated in Fig.~\ref{fig07},
with $\partial A_i=\partial A'_i$. We then have
$$X(\mc{P})=V'_0\cup_{A'_0}T_0\times I\cup_{A'_1}V'_1,$$
where $V'_i=V_i\setminus\intr N(\mc{P})\subset V_i$ is a solid
torus. Now, for $i=0,1$, the circle $\wh{c}_i$ is isotopic to a
core of the annulus $A_i\subset\partial X(\mc{P})$; moreover,
if $\mu(\wh{c}_i)\geq 2$ in $X(\mc{P})$ then the annulus $A'_i$
is a companion annulus for a core of $A_i$, so we can take
$V(\wh{c}_i)=V'_i$, while if $\mu(\wh{c}_i)=1$ in $X(\mc{P})$
then the annuli $A_i$ and $A'_i$ are isotopic in $V'_i\subset
X(\mc{P})$. It is not hard to see now that the manifold
$M(\mc{P},\wh{c}_0,\wh{c}_1)$ may always be isotoped in $S^3$
onto the manifold $T_0\times I$ in such a way that $\partial
D\subset\partial M(\mc{P},\wh{c}_0,\wh{c}_1)$ isotopes to a
circle in $\partial (T_0\times I)$ isotopic to $(\partial
T_0)\times\{0\}$, and that such an isotopy can be extended to
an isotopy that maps $M(\mc{P},\wh{c}_0,\wh{c}_1)\cup N(D)$
onto $T\times I$.
\epf

We are now ready to give necessary and sufficient conditions
for a Seifert Klein bottle $\mc{P}$ to be of the form
$P(c_0^{\varepsilon_0}, c_1^{\varepsilon_1},R)$.

\begin{lem}\label{lem2}
Let $K\subset S^3$ be a knot, $P$ any Seifert Klein bottle in
$X_K$ with atoroidal exterior $X(P)$, $c_0,c_1$ any two
disjoint centers of $P$, $\wh{c}_0,\wh{c}_1\subset T_P$ any two
twisted lifts of $c_0,c_1$, respectively, and $D$ the unique
waist disk of $N(P)$ that separates $\wh{c}_0$ and $\wh{c}_1$.
Then, $K$ admits a $\{P,c_0,c_1\}$-structure with induced
twisted lifts $\wh{c}_0,\wh{c}_1 $ iff the following conditions
hold:
\ben
\ita
the manifold $M(P,\wh{c}_0,\wh{c}_1)$ is a genus two
handlebody, and

\itb
$M(P,\wh{c}_0,\wh{c}_1)\cup N(D)\approx\mc{T}\times I$, where
$\mc{T}$ is a closed torus.
\een
\end{lem}

\bpf
If such a $\{c_0,c_1\}$-structure exists, that (a) and (b) hold
follows from Lemma~\ref{mp}.

Conversely, suppose (a) and (b) hold. By (b),
$M(P,\wh{c}_0,\wh{c}_1)\cup N(D)= T\times I$ for some closed
torus $T=T\times \{0\}\subset S^3$; the identity in
(\ref{eqn3}) implies that $T$ is unknotted in $S^3$ and hence
that $M(P,\wh{c}_0,\wh{c}_1)\cup N(D)$ is the exterior of the
Hopf link in $S^3$. We will assume that $\wh{c}_i\subset
T\times\{i\}$ for $i=0,1$, and that $T\times I$ lies inside a
slightly larger product of the form $T\times[-\delta,1+\delta]$
for some sufficiently small $\delta>0$.

Let $N(D)$ be a small regular neighborhood of $D$ in $N(P)$
which is disjoint from $c_i$ and $\wh{c}_i$ for $i=0,1$, and
such that $R=N(D)\cap P$ is a rectangle properly embedded in
$N(D)$ whose core $\alpha\subset R$ is also a cocore of the
2-handle $N(D)$.

Recall from the construction of $M(P,\wh{c}_0,\wh{c}_1)$ that
$\clo(N({P})\setminus N(D))=W_0\sqcup W_1$, where $W_0,W_1$ are
solid tori in $N(P)$ with $\wh{c}_0\subset
\partial W_0$ and $\wh{c}_1\subset\partial W_1$, so that
for $i=0,1$, we have
$$A_i=W_i\cap\left(M({P},\wh{c}_0,\wh{c}_1)\cup N(D)\right)=W_i\cap
(T\times I)=(\partial W_i)\cap \partial(T\times I)$$ is an
annulus with core $\wh{c}\,'_i\subset \intr A_i$ isotopic to
$\wh{c}_i$ in $\partial W_i$ and $\partial (T\times I)$. We may
further assume, after an isotopy of $W_0\cup(T\times I)\cup
W_1$ which leaves $T\times I$ fixed, that
$W_0=A_0\times[-\delta,0],$ and $W_1=A_1\times[1,1+\delta]$ in
$T\times[-\delta,1+\delta]$.

Consider now the Moebius bands $B_i=P\cap W_i$, $i=0,1$; the
rectangle $R$ has one boundary side on $\partial B_0\cap
T\times\{0\}$ and its opposite side on $\partial B_1\cap
T\times\{1\}$, with $P=B_0\cup R\cup B_1$. Since $\wh{c}_0$ is
a twisted lift of $c_0$ in $N(P)$, giving one half-twist to
$B_0$ in $W_0$ away from $R\subset N(D)$ produces an annulus
properly embedded in $W_0$ with the same boundary slope as the
core $\wh{c}\,'_0$ of $A_0$; this annulus can be isotoped
within $W_0$ onto the annulus $\wh{c}\,'_0\times [-\delta,0]$,
with $c_0\subset B_0$ corresponding to the circle
$\wh{c}\,'_0\times \{-\delta/2\}$. In a similar way, one half
twist on $B_1\subset W_1$ (away from $R$) may be assumed to
produce the annulus $\wh{c}\,'_1\times [1,1+\delta]$, with
$c_1\subset B_1$ corresponding to the circle $\wh{c}\,'_1\times
\{1+\delta/2\}$. Finally, the rectangle $R\subset T\times I$
may be slightly isotoped within $T\times I$ so that its side
$R\cap T\times\{0\}$ lies on $\wh{c}\,'_0$ and its side $R\cap
T\times\{1\}$ on $\wh{c}\,'_1$, in such a way that $P$ may be
recovered from the pair of pants $Q=\wh{c}\,'_0\times
[-\delta,0]\cup R\cup\wh{c}\,'_1\times [1,1+\delta]\subset
N(P)$ by performing the corresponding reverse one half twists
on the annuli $\wh{c}\,'_0\times [-\delta,0]\subset W_0$ and
$\wh{c}\,'_1\times [1,1+\delta]\subset W_1$.

Now, the endpoints $q_0=\alpha\cap T\times\{0\}, q_1=\alpha\cap
T\times\{1\}$ of $\alpha$ lie on $\wh{c}\,'_0,\wh{c}\,'_1$,
respectively. Therefore, if $\alpha'$ is the arc
$\{q_0\}\times[-\delta/2,0]\cup\alpha\cup\{q_1\}\times[1,1+\delta/2]\subset
Q$ then $\alpha'$ is an arc properly embedded in the product
$T\times[-\delta/2,1+\delta/2]$ with endpoints on $c_0\sqcup
c_1$ such that $\clo(T\times[-\delta/2,1+\delta/2]\setminus
N(\alpha'))$ is homeomorphic to $\clo(T\times[0,1]\setminus
N(\alpha))$.

By (a), the manifold $\clo(T\times[0,1]\setminus
N(\alpha))\approx \clo(T\times[0,1]\setminus
N(D))=M({P},\wh{c}_0,\wh{c}_1)$ is a genus two handlebody.
Since $T\times[-\delta/2,1+\delta/2]$ and $T\times[0,1]$ are
isotopic in $S^3$, and $T\times[0,1]\subset S^3$ is the
exterior of the Hopf link, it follows that $\alpha'$ is a
tunnel for the Hopf link exterior
$T\times[-\delta/2,1+\delta/2]\subset S^3$. As the Hopf link is
a 2-braid link, by \cite[Theorem 2.1]{reid1} such a tunnel arc
$\alpha'$ is unique up to isotopy, hence $\alpha'$ can be
isotoped into a fiber $\{*\}\times [-\delta/2,1+\delta/2]$ of
$T\times[-\delta/2,1+\delta/2]$. Therefore $P$ is isotopic to a
surface of the form $P(c_0^{\varepsilon_0},
c_1^{\varepsilon_1},R)$ and so $K$ admits a
$\{P,c_0,c_1\}$-structure.
\epf

With the aid of Lemma~\ref{comp2}, the conditions in
Lemma~\ref{lem2} can now be expressed entirely in simple
algebraic terms whenever the Seifert Klein bottle surface $P$
is unknotted; this is the main content of the next result.

\begin{lem}\label{lem3}
Suppose $P$ is an unknotted Seifert Klein bottle for a knot
$K\subset S^3$. Let $c_0,c_1$ be disjoint centers of $P$ with
corresponding twisted lifts $\wh{c}_0,\wh{c}_1\subset\partial
X(P)$ and $D$ the unique waist disk of $N(P)$ that separates
$\wh{c}_0$ and $\wh{c}_1$.
\ben
\ita
The knot $K$ admits a $\{P,c_0,c_1\}$-structure with induced
twisted lifts $\wh{c}_0,\wh{c}_1$ iff one of the following set
of conditions holds:
\ben
\item[(i)]
$\mu(\,\wh{c}_0)=1=\mu(\,\wh{c}_1)$ and $[\partial
D]\equiv[u,v]$ in $\pi_1(X(P))$ for some (and hence any) basis
$\{u,v\}$ of $\pi_1(X(P))$;

\item[(ii)]
for some $\{i,j\}=\{0,1\}$: $\mu(\,\wh{c}_i)=1$,
$\mu(\,\wh{c}_j)\geq 2$, and $\partial
D\equiv[u,v^{\mu(\,\wh{c}_j)}]$ in $\pi_1(X(P))$ for some (and
hence any) basis $\{u,v\}$ of $\pi_1(X(P))$ such that
$[\,\wh{c}_j]\equiv v^{\mu(\,\wh{c}_j)}$;

\item[(iii)]
$\mu(\,\wh{c}_0)\geq 2$, $\mu(\,\wh{c}_1)\geq 2$, and $\partial
D\equiv[u^{\mu(\,\wh{c}_0)},v^{\,\mu(\wh{c}_1)}]$ in
$\pi_1(X(P))$ for some (and hence any) basis $\{u,v\}$ of
$\pi_1(X(P))$ such that $[\,\wh{c}_0]\equiv
u^{\mu(\,\wh{c}_0)}$ and $[\,\wh{c}_1]\equiv
v^{\mu(\,\wh{c}_1)}$.
\een

\itb
Suppose $K$ admits a $\{P,c_0,c_1\}$-structure with induced
twisted lifts $\wh{c}_0,\wh{c}_1$, and let $\{i,j\}=\{0,1\}$.
If $\mu(\,\wh{c}_i)=1$ in $X(P)$ then, in $\pi_1(X(P))$, either
$[\,\wh{c}_i],[\,\wh{c}_j]$ are both primitive or
$[\,\wh{c}_j]$ is a power of a primitive.
\een
\end{lem}

\bpf
As $X(P)$ is a genus two handlebody, part (a) follows by a
direct application of Lemma~\ref{comp2} to Lemma~\ref{lem2}; in
particular, observe that, by Lemma~\ref{comp2}, a basis
$\{u,v\}$ of $\pi_1(X(P))$ satisfying the condition
$[\,\wh{c}_j]\equiv v^{\mu(\,\wh{c}_j)}$ in (ii) or the
conditions $[\,\wh{c}_0]\equiv u^{\mu(\,\wh{c}_0)}$ and
$[\,\wh{c}_1]\equiv v^{\mu(\,\wh{c}_1)}$ in (iii) always
exists.

For part (b) we have that $P$ is a surface of the form
$P(c_0^{\varepsilon_0}, c_1^{\varepsilon_1},R)$ with induced
twisted lifts $\wh{c}_0,\wh{c}_1$; we assume for definiteness
that $\mu(\wh{c}_0)=1$ in $X(P)$. If $\mu(\wh{c}_1)=1$ too then
$X(P)= M(P,\wh{c}_0,\wh{c}_1)$; since, by Lemma~\ref{mp},
$M(P,\wh{c}_0,\wh{c}_1)=T_0\times I$ for some once punctured
torus $T_0\subset S^3$ such that $\wh{c}_0\subset T_0\times
\{0\}$ and $\wh{c}_1\subset T_0\times \{1\}$, it follows that
$[\wh{c}_0],[\wh{c}_1]$ are primitive in
$\pi_1(X(P))=\pi_1(T_0\times I)=\pi_1(T_0)$. Otherwise,
$\mu(\wh{c}_1)\geq 2$ and so $[\wh{c}_1]$ is a power of a
primitive in $\pi_1(X(P))$ by Lemma~\ref{comp}(b).
\epf

\subsection{A tunnel number one criterion}
The next result gives a condition for a knot with an unknotted
Seifert Klein bottle to have tunnel number one.

\begin{lem}\label{tunnel}
Let $K\subset S^3$ be a nontrivial knot with an unknotted
Seifert Klein bottle $P$. Then $K$ has tunnel number one if
some center circle of $P$ has a lift to $T_P$ which is
primitive in $X(P)$.
\end{lem}

\bpf
Since $P$ is unknotted, it follows from (\ref{basic}) that
$N(P)\cup X(P)$ is a genus two Heegaard decomposition of $S^3$.
Let $c_0$ be any center circle of $P$, and let $c_1$ be the
only other center circle of $P$ which is disjoint from $c_0$.
Then $N(P)=W_0\cup_D W_1$, where $W_0,W_1\subset N(P)$ are
solid tori and $D$ is the (unique) waist disk $D$ of $N(P)$
which separates $c_0,c_1$. Finally, let $B(c_0)$ be the Moebius
band in $N(P)$ constructed in Section \ref{kbottles}, so that
$\partial B(c_0)\subset T_P$ is the lift of $c$, and let
$V(c_0)$ be a regular neighborhood of $B(c_0)$ in $N(P)$
disjoint from $\partial P$; $V(c)$ is a solid torus.

If $\partial B(c_0)$ is primitive in $X(P)$ then $H_1=X(P)\cup
V(c_0)$ is a genus two handlebody, while clearly
$H_2=\clo\left(N(P)\setminus V(c_0)\right)\subset S^3$ is also
a genus two handlebody; thus $H_1\cup H_2$ is a Heegaard
decomposition for $S^3$. From the representation of
$N(P),D,\partial P$ shown in Fig.~\ref{fig13} it is not hard to
see that $[\partial P]$ is primitive in $\pi_1(H_2)$, which
implies that there is a properly embedded disk $D_0$ in $H_2$
which intersects $\partial P$ transversely in one point.
Therefore the knot $\partial P$, and hence $K$, has tunnel
number one.
\epf

\subsection{Hyperbolicity}

\begin{lem}\label{lem4}
Let $k\subset S^3$ be  a knot which is a (possibly trivial)
torus knot or a 2-bridge knot. Let $P$ be a Seifert Klein
bottle for $k$ with meridian circle $m$, and let $T_P$ be the
frontier of $N(P)$ in $X_k$.
\ben
\ita If $k$ is a trivial or 2-bridge knot then $m$ is a trivial knot in
$S^3$.

\itb
If $k$ is nontrivial then $P$ is unknotted; moreover, if $k$ is
a torus knot then
\ben
\item[(i)]
$T_P$ is not $\pi_1$-injective in $X_k$,

\item[(ii)]
$T_P$ is boundary compressible in $X_k$ iff $k=T(3,5), \
T(3,7)$, or $T(2,n)$ for some odd integer $n\neq\pm 1$, in
which case
$$ m=\begin{cases}
\text{trivial knot,} & k\neq T(3,7)\\
T(2,3), & k= T(3,7).
\end{cases}
$$
\een
\een
\end{lem}

\bpf
If $k$ is a 2-bridge knot then it follows from \cite[Theorem
1.2(a)]{valdez8} that $k$ is a plumbing of an annulus $A$ and a
Moebius band $B$, which are unknotted and unlinked in $S^3$.
Moreover, $k$ bounds a unique Seifert Klein bottle $P$ by
\cite{thurs4}, namely the one obtained as the plumbing of $A$ and
$B$, so $P$ is unknotted. This proves part (a) for $k$ a
2-bridge knot; since the meridian circle $m$ of $P$ is the core
of $A$, it also follows that $m$ is a trivial knot in $S^3$.

The only nontrivial torus knots that bound a Moebius band are
those of the form $T(2,n)$ for some odd integer $n\neq \pm 1$.
Hence, any nontrivial torus knot $k$ which is not of the form
$T(2,n)$ and bounds a Seifert Klein bottle has crosscap number
two, so it follows from \cite[Corollary 1.6]{valdez7} that $P$
is unique in $X_k$ up to isotopy and that $P$ is unknotted and
not $\pi_1$-injective in $X_k$; moreover, by the argument used
in the proof of \cite[Corollary 1.6]{valdez7}, $T_P$ is
boundary compressible in $X_k$ iff $k$ is a torus knot of the
form $T(5,3),T(3,7)$. Projections of the knots $T(3,5),T(3,7)$
spanning their unique Seifert Klein bottle $P$ are shown in
Fig.~\ref{fig08}; it can be seen that indeed the meridian
circle $m$ of $P$ is the trivial knot and the trefoil $T(2,3)$,
respectively ($T(3,5)$ can also be drawn as the pretzel
$(-2,3,5)$, again showing the meridian circle of its Seifert
Klein bottle is trivial).

\begin{figure}
\Figw{5in}{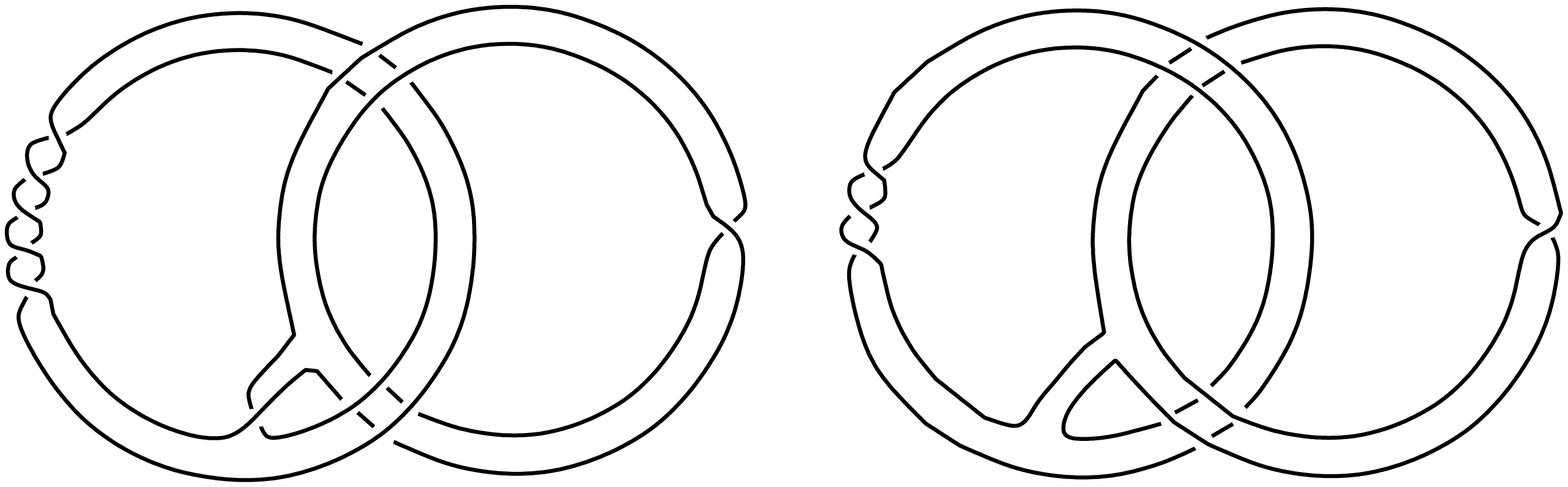}{The torus knots $T(3,7)$ and $T(3,5)$.}{fig08}
\end{figure}

We now deal with the cases where $k$ is either a nontrivial
torus knot of the form $T(2,n)$, or $k$ is the trivial knot
with a Seifert Klein bottle whose boundary slope, relative to a
standard meridian-longitude pair $\mu,\lambda\subset\partial
X_k$, is not $0/1$. These cases have in common that any Seifert
Klein bottle $P$ in $X_k$, and hence the associated surface
$T_P$, boundary compress in $X_k$: this follows from
\cite[Lemma 4.2]{valdez6} if $k=T(2,n)$, while if $k$ is
trivial then $X_k$ is a solid torus and hence the surface $P$
must boundary compress in $X_k$. In either case, $P$ boundary
compresses into a Moebius band $B\subset X_k$ such that
$\Delta(\partial P,\partial B)=2$. If $k=T(2,n)$ then $B$ is
unique and so, for a fixed boundary slope, $P$ is obtained by
adding a band (rectangle) in $\partial X_k$ to the Moebius band
$B$, and hence $P$ is unique up to isotopy in $X_k$; if $k$ is
the trivial knot we will see below that though $B$ may not be
unique, the surfaces $P$ are unique relative to their boundary
slopes.

Using a standard meridian-longitude pair $\mu,\lambda$ in
$\partial X_k$, we may assume that $\partial P=a\mu+\lambda$
and
$$\partial B=\begin{cases}
2\mu+b\lambda \ (b=\text{odd}), & k=\text{trivial}\\
2n\mu+\lambda, & k=T(2,n).
\end{cases}$$
Since $\Delta(\partial P,\partial B)=2$, we must then have
$$
a=\begin{cases}
\pm\ 4\text{ and }b=\mp 1, & k=\text{trivial}\\
2n\pm 2, & k=T(2,n).
\end{cases}
$$

The knot $k=T(2,n)$ can be represented as the pretzel knots
$(n-2,1,-2)$ or $(n+2,-1,2)$, and in these projections $k$
spans a Seifert Klein bottle with boundary slope $2n-2,2n+2$,
respectively, all having meridian circle a trivial knot. Since
the Seifert Klein bottle for $k$ is unique for each of these
slopes, it follows that the meridian circle of any Seifert
Klein bottle bounded by $k=T(2,n)$ is a trivial knot.

If $k$ is the trivial knot then it follows from the above
computation that $k$ bounds two Seifert Klein bottles with
boundary slopes distinct from $0/1$: $P_1$ with boundary slope
$+4$, which boundary compresses to the unique Moebius band
$B_1\subset X_k$ such that $\partial B_1=2\mu-\lambda$, and
$P_2$ with boundary slope $-4$, which boundary compresses to
the unique Moebius band $B_2\subset X_k$ such that $\partial
B_2=2\mu+\lambda$. Therefore $P_1$ and $P_2$ are unique
relative to their boundary slope; the trivial knot $k$ can be
represented as the pretzel knots $(0,1,1)$ or $(0,-1,-1)$ and
in these projections $k$ spans the unique Seifert Klein bottles
with boundary slopes $4,-4$, respectively, both of which have
the trivial knots as meridian circle.

Therefore part (b) holds, and for part (a) only the case when
$k$ is trivial and $P$ has boundary slope $0/1$ remains. But in
this last case, if $D$ is the meridian disk of the solid torus
$X_k$ then $\partial D$ and $\partial P$ both have the same
boundary slope and so $P$ and $D$ can be isotoped in $X_k$ so
that $\partial D$ and $\partial P$ are disjoint, and hence so
that $P\cap D$ intersect transversely and minimally in circles
only. Necessarily, $P\cap D$ is nonempty and consists of
orientation preserving circles in $P$, ie of meridians or
longitudes of $P$ (cf Section \ref{kbottles}). As an innermost
disk $D'$ of $P\cap D\subset D$ compresses $P$ along $\partial
D'$, the circle $\partial D$ must be a meridian of $P$, for
otherwise $\partial D$ would be a longitude of $P$ that bounds
some Moebius band $B\subset P$, and $B\cup_{\partial} D'$ would
be a closed projective plane in $S^3$, an impossibility. Hence
the meridian circle of $P$ is a trivial knot in $S^3$.
\epf

The above lemma can be used to give conditions under which a
knot bounding a Seifert Klein bottle $P$ is hyperbolic in terms
of the lifts of the meridians and the longitudes of $P$.

\begin{lem}\label{hyper1}
Let $K$ be a nontrivial knot with Seifert Klein bottle $P$ such
that $X(P)$ is atoroidal. Let $m_1,m_2$ be the lifts of the
meridian circle $m$ of $P$ and $\mc{L}$ the collection of lifts
of longitudes of $P$. If, in $X(P)$, $\mu(\ell')=1$ for each
$\ell'\in\mc{L}$ and either $\mu(m_1)=1$ or $\mu(m_2)=1$, then
$X_K$ is atoroidal, and in such case,
\ben
\ita
if $P$ is $\pi_1$-injective in $X_K$ then $K$ is hyperbolic,

\itb if the frontier $T_P$ of $N(P)$ is boundary compressible in
$X_K$ and the meridian circle $m$ of $P$ is a nontrivial knot
in $S^3$ then either $K$ is hyperbolic and not 2-bridge, or
$K=T(3,7)$ and $m=T(2,3)$.
\een
\end{lem}

\bpf
By Lemma~\ref{comp}(b), no element of $\mc{L}$, nor say $m_2$,
has a companion annulus in $X(P)$.

Suppose $T$ is an essential torus in $X_K$. Since, by
\cite[Lemma 4.2]{valdez6}, the fact that $K$ is a nontrivial knot
implies that the surface $P$ is incompressible in $X_K$, and
both $N(P)$ and $X(P)$ are atoroidal, we may assume $T$
intersects $P$ transversely and minimally with $P\cap T$ a
nonempty collection of circles which are nontrivial in both $P$
and $T$. Necessarily, $P\cap T\subset T$ is a family of
mutually parallel nontrivial circles, so we may assume that
$T\cap N(P)$ is a collection of annuli which are fibered under
the $I$-bundle structure of $N(P)=P\wt{\times}I$, ie $T\cap
N(P)=(T\cap P)\wt{\times}I\subset N(P)$, and that each
component of $T\cap X(P)$ is an annulus. Moreover, each circle
in $P\cap T\subset P$ is either a meridian circle, a longitude
circle, or a circle parallel to $\partial P$. It is not hard to
see that if all circles $P\cap T\subset P$ are parallel to
$\partial P$ then $T$ must be parallel to $\partial X_K$, which
is not the case. Therefore there is a component $A$ of $T\cap
X(P)$ such that at least one of its boundary components
$\partial _1 A,\partial _2 A$ is not parallel to $\partial P$
in $P$. By \cite[Lemma 2.3]{valdez7} we then have that, in
$\partial X(P)$, the circles $\partial _1 A,\partial _2 A$ are
either both parallel to $m_1$, both parallel to $m_2$, or both
parallel to some element of $\mc{L}$. Given our hypothesis on
$\mc{L}$ and $m_1,m_2$, since $|P\cap T|$ is minimal, this
implies necessarily that $A$ is a companion annulus in $X(P)$
for $m_1$. By minimality of $|P\cap T|$, there are two annular
components $A',A''$ of $T\cap N(P)$ with, say, $\partial_1
A=\partial_1 A', \partial_2 A=\partial_1 A''$, and with
$\partial_2 A,\partial_2 A''\subset\partial N(P)=\partial X(P)$
parallel to $m_2$. Applying \cite[Lemma 2.3]{valdez7} again, it
follows that some component of $T\cap X(P)$ must be a companion
annulus of $m_2$ in $X(P)$, which by hypothesis is not the
case.

Therefore $X_K$ is atoroidal and so by \cite{thurs2} the knot
$K$ is either a hyperbolic or torus knot, hence part (a) holds
by Lemma~\ref{lem4}(b)(i). Suppose now $T_P$ is boundary
compressible in $X_K$. If $m$ is a nontrivial knot then $K$ is
not a 2-bridge knot by Lemma~\ref{lem4}(a), and if $K$ is a
torus knot then by Lemma~\ref{lem4}(b)(ii)  $K$ must be the
$T(3,7)$ torus knot and $m=T(2,3)$; hence (b) holds.
\epf

\subsection{Meridians, centers, longitudes, and their lifts}

Our first result gives a criterion to identify the lifts of a
meridian or longitude circle of a once punctured Klein
bottle$P$; along with Lemma~\ref{comp}(c), this result will
provide a simple algebraic way of identifying the lifts of the
meridian.

\begin{lem}\label{lem7}
Let $P$ be a once punctured Klein bottle and $H=P\wt{\times}I$.
If $A$ is an incompressible annulus properly embedded in $H$
with $\partial A\subset \partial H\setminus\partial P$, and $A$
is not parallel into $\partial H$, then $A$ can be isotoped in
$H$ so that $A\cap P$ is either one meridian (if $A$ is
nonseparating) or one longitude (if $A$ is separating) circle
of $P$; in particular, $\partial A$ are the lifts of the circle
$A\cap P$ to $\partial H$.
\end{lem}

\bpf
Let $A_P$ be an annular regular neighborhood of $\partial P$ in
$\partial H$, $T$ be the twice punctured torus $\partial
H\setminus\intr A_P$, and $M$ be the manifold obtained by
cutting $H$ along $P$, so that $M=T\times I$ with $T$
corresponding to $T\times 0$.

Isotope $A$ and $P$ in $H$ so as to intersect transversely and
minimally. If $A\cap P=\emptyset$ then $A$ lies in $M=T\times
I$ with $\partial A\subset T\times 0$; but then $A$ is parallel
into $T\times 0$ in $M$ (cf \cite[Corollary 3.2]{wald1}), and
hence $A$ is parallel into $\partial H$ in $H$, contradicting
our hypothesis; thus $A\cap P\neq\emptyset$. Since $A$ is
orientable, $A\cap P$ consists of a collection of circles which
preserve orientation in $P$, ie, of longitudes of $P$ only or
meridians of $P$ only. If $\alpha,\beta$ are distinct
components of $A\cap P$ which cobound an annulus $A'\subset A$
with $P\cap\intr A'=\emptyset$, then the annulus $A'\cap M$ has
both of its boundary components in $T\times 1$, and hence $A'$
is parallel into $T\times 1$ in $T\times I$. It follows that
$A'$ is parallel into $P$ in $H$, contradicting the fact that
$A\cap P$ is minimal. Therefore, $A\cap P$ consists of a single
circle $\gamma$, a meridian or longitude of $P$, and so
$\partial A$ are the lifts of $\gamma$ to $\partial H$.
\epf

Now let $m$ be the meridian circle of $P$, and let $c,\ell$ be
a center and a longitude circle of $P$, respectively. If
$c(n),\ell(n)\subset P$ denote the circles obtained by
Dehn-twisting $n$-times the circles $c,\ell$ along $m$, it
follows from \cite[Lemma 3.1]{valdez7} that the collections
$\Set{c(n)}{n\in\mathbb{Z}}$ and
$\Set{\ell(n)}{n\in\mathbb{Z}}$ consist of all center and
longitude circles of $P$ up to isotopy.

This fact generalizes into the following result, which
describes the construction of the twisted lifts of all centers
of $P$ and of all lifts of longitudes; its proof follows easily
from the fact that Dehn-twisting $P\wt{\times}I$ along the
fibered annulus $A(m)\subset P\wt{\times}I$ defined in Section
\ref{kbottles} is an automorphism of $H$ that fixes $m$ and
maps $P\subset H$ into itself.

\begin{lem}\label{311}
Let $P$ be a once punctured Klein with meridian circle $m\in P$
and fibered annulus $A(m)\subset H=P\wt{\times}I$. Let
$c_0,c_1$ be a pair of disjoint center circles of $P$ with
$c_0',c_0''$ and $c_1',c_1''$ their twisted lifts to $\partial
H$, respectively, and let $\ell'\subset\partial H$ be the lift
of any longitude of $P$.

Let $c_0(n),c_1(n)\subset P$ and
$c'_0(n),c''_0(n),c'_1(n),c''_1(n),\ell'(n)\subset\partial H$
be the circles obtained from
$c_0,c_1,c'_0,c''_0,c'_1,c''_1,\ell'$ after Dehn twisting $H$
$n$-times along $A(m)$. Then,

\ben
\ita
the collection $\Set{(c_0(n),c_1(n))}{n\in\mathbb{Z}}$ consists
of all pairs of disjoint centers of $P$, and the twisted lifts
of each $c_i(n)$ are the circles $c_i'(n),c_i''(n)$;

\itb
the collection $\Set{\ell'(n)}{n\in\mathbb{Z}}$ consists of all
lifts of longitudes of $P$.
\hfill\qed
\een
\end{lem}

\section{The knots \texorpdfstring{$K(p,q)$}{K(p,q)}}\label{sec4}

We begin this section by providing a detailed construction of
the family of knots $K(p,q)$ given in Section
\ref{introduction}. Let $H$ be a genus two handlebody
standardly embedded in $S^3$ with a complete system of meridian
disks $D_0,D_1$, as shown in Fig.~\ref{fig10}, and let
$c_0',c_0''$ and $c_1',c_1''$ be the four circles embedded in
the boundary $\partial H$ shown in the same figure; notice that
any pair formed by one circle in $c_0',c_0''$ and one in
$c_1',c_1''$ give rise to a basis of $\pi_1(H)$.
\begin{figure}
\begin{center}
\psfrag{c0}{$c'_{0}$}
\psfrag{c1}{$c''_0$}
\psfrag{c2}{$c'_1$}
\psfrag{c3}{$c''_1$}
\psfrag{D0}{$D_0$}
\psfrag{D1}{$D_1$}
\Figw{4.5in}{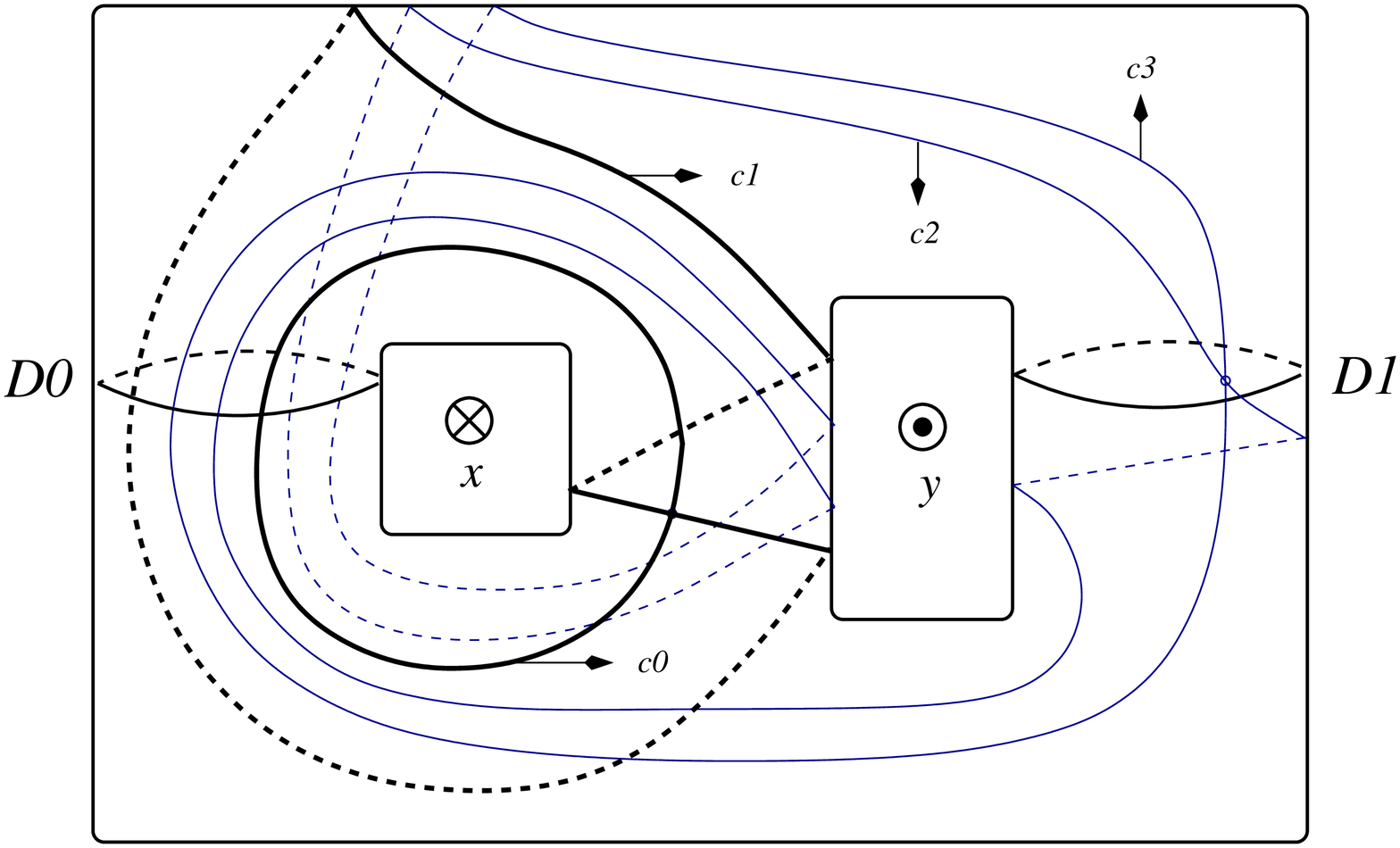}{}{fig10}
\end{center}
\end{figure}
The homological sums $\ell'_0$ and $\ell'_1$ of the pairs of
curves $c_0',c_0''$ and $c_1',c_1''$, respectively, indicated
in Fig.~\ref{figmoebius}, bound disjoint Moebius bands $B_0,
B_1$ in $H$, respectively.

\begin{figure}
\begin{center}
\psfrag{B0}{$B_0$}
\psfrag{B1}{$B_1$}
\psfrag{l0}{$\ell_0'$}
\psfrag{l1}{$\ell_1'$}
\psfrag{D1}{$D_1$}
\psfrag{D0}{$D_0$}
\Figw{4.5in}{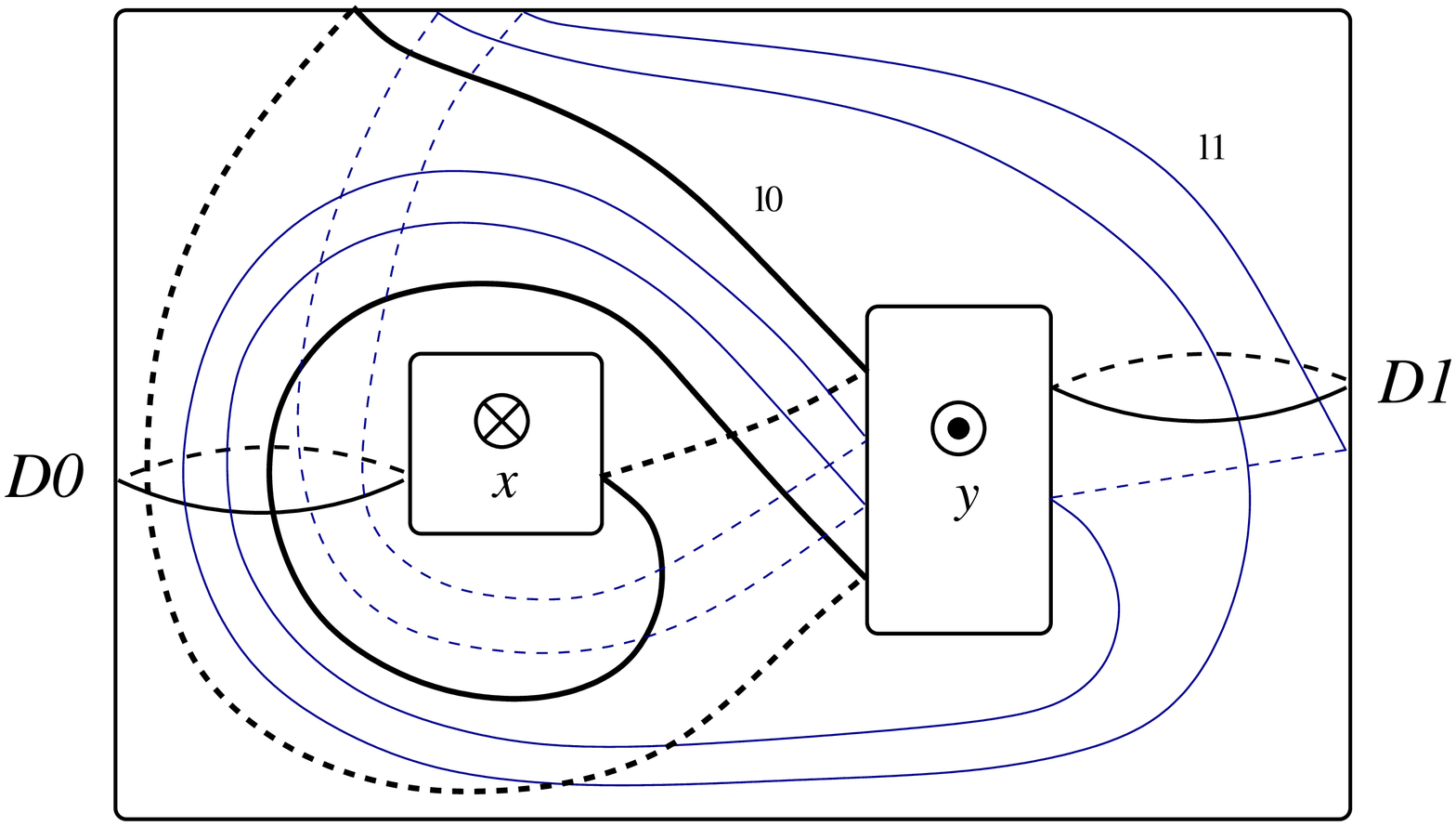}{}{figmoebius}
\end{center}
\end{figure}

\begin{figure}
\begin{center}
\psfrag{D0}{$D_0$}\psfrag{D1}{$D_1$}
\psfrag{H}{$H$}\psfrag{R}{$R$}
\psfrag{Test}{$K(0,0)$}
\Figw{4.5in}{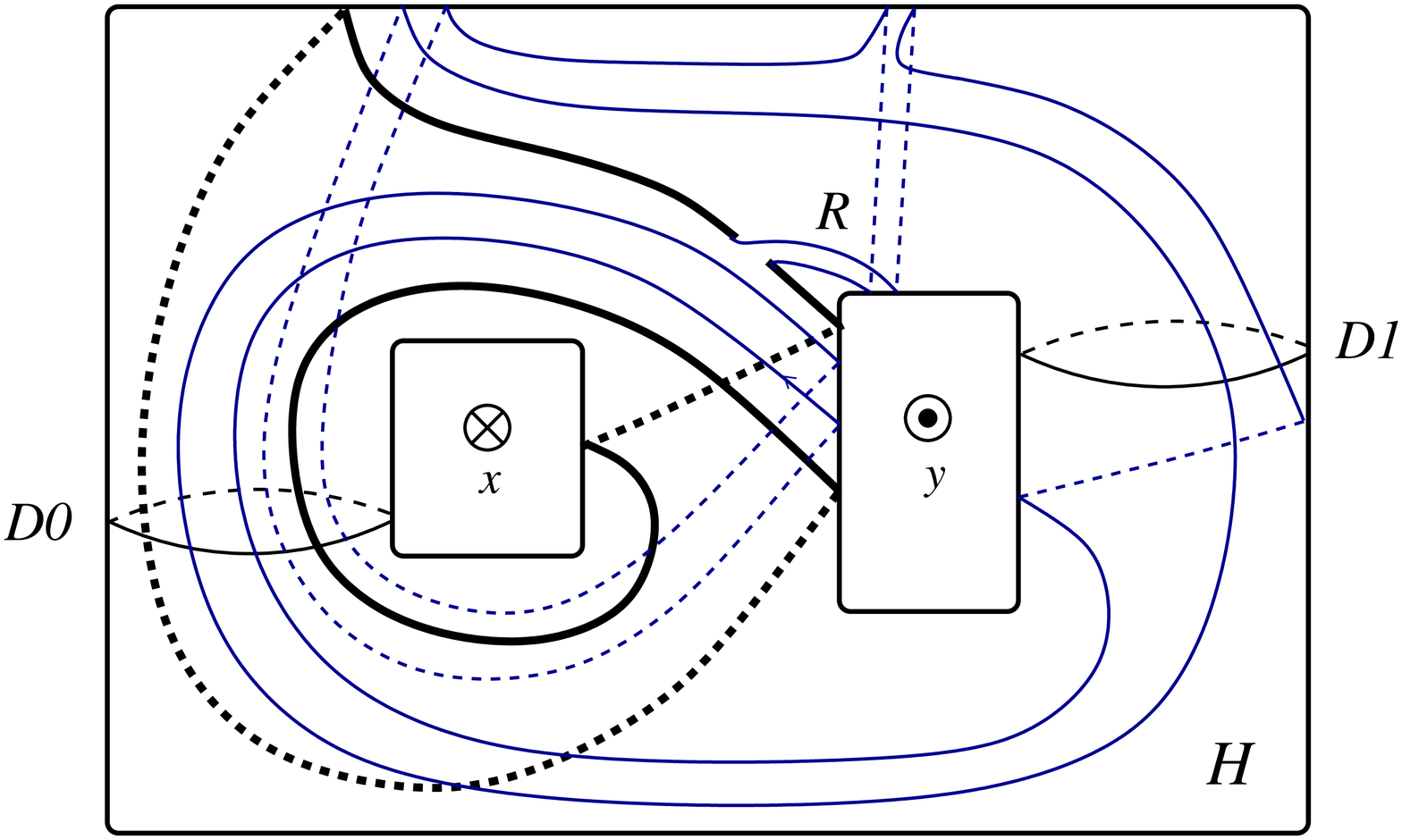}{The Seifert Klein bottle
bounded by $K(0,0)$ in $H$.}{fig11a}
\end{center}
\end{figure}

A waist disk $D_2$ of $H$ separating $B_1$ and $B_2$ can be
constructed from the meridian disk $D_0$ of $H$ in
Fig.~\ref{figwaistdisk}, which is disjoint from $B_1$ and
intersects $B_0$ in a single essential arc, by taking the
frontier of a regular neighborhood of $D_0\cup c'_1$ (or of
$D_0\cup c''_1$) in $H$; such a waist disk of $H$ is unique up
to isotopy.

\begin{figure}
\begin{center}
\psfrag{D}{$D_2$}\psfrag{D0}{$D_0$}\psfrag{H}{$H$}
\Figw{4in}{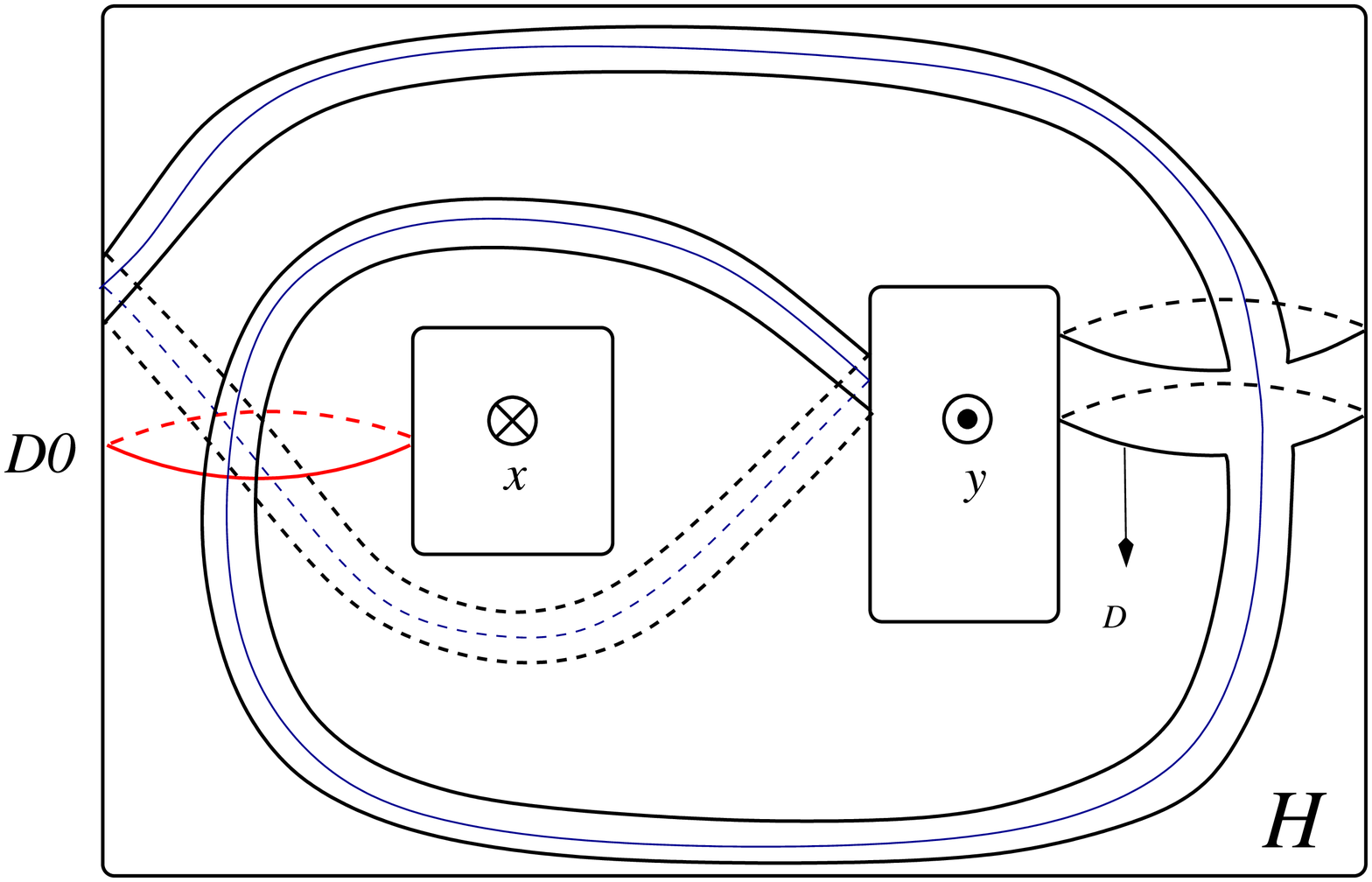}{}{figwaistdisk}
\end{center}
\end{figure}

We now connect the Moebius bands $B_0,B_1$ with the rectangle
$R$ shown in Fig.~\ref{fig11a} and produce a properly embedded
Seifert Klein bottle $P(0,0)$ in $H$ whose boundary $\partial
P(0,0)$ is isotopic to the knot $K(0,0)$ shown in
Fig.~\ref{fig11b}; notice the cores $c_0,c_1$ of $B_0,B_1$,
respectively, are disjoint centers of $P(0,0)$. The fact that
$R$ intersects $D_2$ transversely in a single arc implies that
$H=P(0,0)\wt{\times} I$, whence $H'=X(P(0,0))$, by Section
\ref{twisted} that the twisted lifts of the centers $c_0,c_1$
are the pairs $c_0',c_0''$ and $c_1',c_1''$, respectively, and
by Section \ref{kbottles} that $\ell_0'=\partial B_0$ and
$\ell_1'=\partial B_1$ are the lifts of the longitudes of $P$
corresponding to the centers $c_0,c_1$, respectively. In
particular, $P(0,0)$ is unknotted.

Finally, consider the circles $m_0,m_1$ shown in
Fig.~\ref{figmer}, which are disjoint from $K(0,0)$. Relative
to the base of $\pi_1(H)$ dual to the meridian disks $D_0,D_1$,
$m_0,m_1$ give rise to conjugate primitive words and so
$m_0,m_1$ cobound a nonseparating annulus in $A(0,0)$ in $H$ by
Lemma~\ref{comp}(c). By Lemma~\ref{lem7}, it follows that
$m_0,m_1$ are the lifts of the meridian $m$ of $P(0,0)$, so
$A(0,0)$ can be isotoped so as to intersect $P(0,0)$ in $m$;
thus $A(0,0)$ is the fibered annulus generated by $m$ in
$H=P\wt{\times} I$.

\begin{figure}
\begin{center}
\psfrag{m0}{$m_0$}
\psfrag{m1}{$m_1$}
\psfrag{D0}{$D_0$}
\psfrag{D1}{$D_1$}
\psfrag{H}{$H$}
\Figw{4.5in}{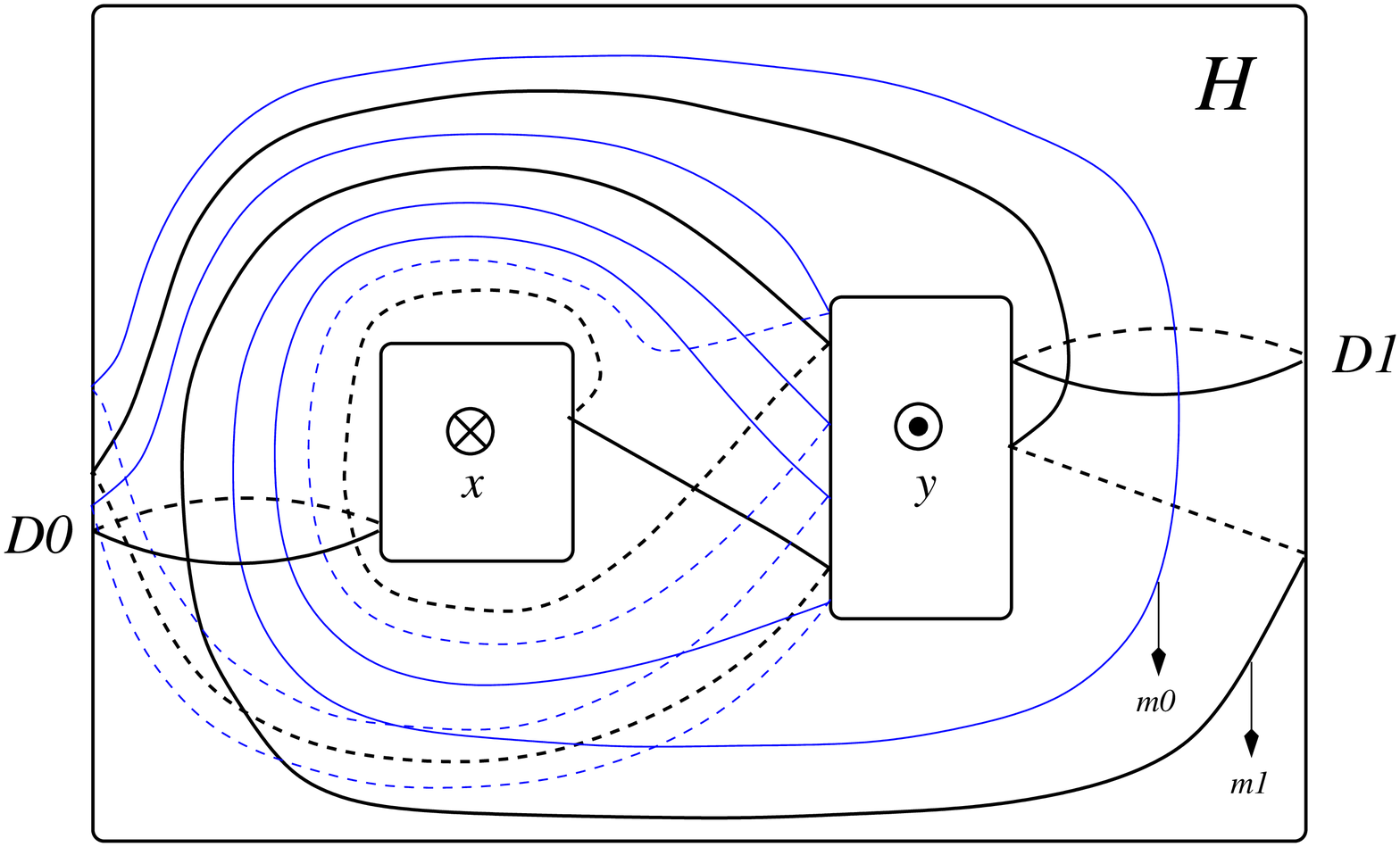}{}{figmer}
\end{center}
\end{figure}

By Lemma~\ref{311}, Dehn twisting $H$ $n$ times along $A(0,0)$
for $n\in\mathbb{Z}$ gives rise to the collections
$\{(c_0'(n),c_0''(n)) \ |  n\in\mathbb{Z}\}$ and
$\Set{(c_1'(n),c_1''(n))}{n\in\mathbb{Z}}$ of twisted lifts of
disjoint pairs of centers of $P$, as well as to the collection
$\Set{\ell'(n)}{n\in\mathbb{Z}}$ of all lifts of longitudes of
$P$ (with $c_0'$ giving rise to $c_0'(n)$, $\ell'_0$ to
$\ell'_0(n)$, etc). At the same time the waist disk $D_2\subset
H$ gives rise to a waist disk $D_2(n)$ separating the twisted
lifts $c_0'(n)\cup c_0''(n)$ and $c_1'(n)\cup c_1''(n)$.

\begin{rem}
Any once punctured Klein bottle with $H$ as regular
neighborhood can be constructed following the procedure
outlined above for $P(0,0)$. That is, once the pairs of circles
$c_0',c_0''$ and $c_1',c_1''$ (such that any pair formed by one
circle in $c_0',c_0''$ and one in $c_1',c_1''$ is a basis for
$\pi_1(H)$), and the disk $D_2$ separating them are given, any
rectangle $R\subset H$ that intersects $D_2$ in one arc may be
used to join the Moebius bands $B_0,B_1$; the latter condition
on $R$ guarantees that $H=P\wt{\times}I$. In general, different
rectangles $R$ will give rise in $H$ to nonisotopic once
punctured Klein bottles whose meridians induce nonisotopic
fibered annuli.
\end{rem}

For any integers $p,q$, the Seifert Klein bottle $P(p,q)$ and
the knot $K(p,q)=\partial P(p,q)$ are then obtained by Dehn
twisting the pair $(P(0,0), K(0,0))$ $p$ and $q$ times along
the disks $D_0,D_1$ of Fig.~\ref{fig11b}, respectively, as
explained in Section \ref{introduction}. We will denote the
pairs of twisted lifts of centers of $P(p,q)$ by
$c_0'(n,p,q),c_0''(n,p,q)$ and $c_1'(n,p,q), c_1''(n,p,q)$, by
$D_2(n,p,q)\subset H$ the waist disk that separates these
pairs, by $\ell'(n,p,q)$ the lifts of the longitudes of
$P(p,q)$, and by $A(p,q)\subset H$ the fibered annulus in
$H=P(p,q)\wt{\times}I$ induced by the meridian circle $m$ of
$P(p,q)$ (with $c'_0(n)$ giving rise to $c'_0(n,p,q)$ after the
Dehn twists along $D_0,D_1$, etc).

By construction we have $H=N(P(p,q))=P(p,q)\wt{\times}I$ and
$H'=X(P(p,q))$, so $P(p,q)$ is unknotted. Let
$\pi_1(H')=\langle x,y\rangle$, where $x,y$ is the base dual to
the complete disk system of $H'$ indicated in Fig.~\ref{fig11b}
(with the disks, and hence the circles $x,y$, oriented by the
arrow head and arrow tail shown in the same figure). It is not
hard to see that the words in $\pi_1(H')=\langle x,y \ | \
-\rangle$ corresponding to the circles constructed above, up to
equivalence, are the ones given in the next lemma; for
convenience, we may denote the elements
$[c_0'(n,p,q)],[c_0''(n,p,q)]$, etc, of $\pi_1(H')$ simply by
$c_0'(n,p,q),c_0''(n,p,q)$, etc.

\begin{lem}\label{words}
In $\pi_1(H')=\langle x,y \ | \ -\rangle$,
\ben
\item[(i)]
$c_0'(n,p,q)\equiv x^p(x^{p+1}yx^py^{q+1}x^py)^n$;

\medskip

\item[(ii)]
$c_0''(n,p,q)\equiv x^pyxy(\ovy x^{\,p}yx^py^{q}x^pyx^p\,\ovy\,
\ovx^{\,p})^n$;

\medskip

\item[(iii)]
$c_1'(n,p,q)\equiv
y\ovx^{\,p}\,\ovy\,\ovx^{\,p}\,\ovy^{\,q}\,(y^{q+1}x^py
x^{p+1}yx^p)^n$;

\medskip

\item[(iv)]
$c_1''(n,p,q)\equiv
\ovy^{\,q}\,\ovx^{\,p}\,\ovy\,\ovx^{\,p}(x^{p}yx^py^{q}x^pyx^p\,\ovy\,
\ovx^{\,p}\,\ovy)^n$;

\medskip

\item[(v)]
$\ell'_0(n,p,q)\equiv (x^{p+1}yx^py^{q+1}x^py)^n x^{p+1}y
(\ovy\, x^{p}yx^py^{q}x^pyx^p\,\ovy\, \ovx^{\,p})^n x^py$;

\medskip

\item[(vi)]
$\ell'_1(n,p,q)\equiv (y^{q+1}x^py x^{p+1}yx^p)^n
y\,\ovx^{\,p}\,\ovy\,\ovx^{\,p}\,\ovy^{\,q}\,
\ovx^{\,p}\,\ovy\,\ovx^{\,p}
(x^{p}yx^py^{q}x^pyx^p\,\ovy\,
\ovx^{\,p}\,\ovy)^n\, \ovy^{\,q}$;

\item[(vii)]
\begin{multline*}
\partial D_2(n,p,q)\equiv (x^py x^{p+1}yx^py^{q+1})^n
(\ovy^{\,q}\,\ovx^{\,p}\,\ovy\,\ovx^{\,p}\,yx^py\,\ovx^{\,p}\,\ovy\,\ovx^{\,p})^n\\
\cdot(yx^py\ovx^{\,p}\,\ovy\,\ovx^{\,p}\,\ovy^{\,q}\,\ovx^{\,p}\,\ovy\,\ovx^{\,p})^n
x^pyx^p(y^qx^pyx^p\,\ovy\,\ovx^{\,p}\,\ovy\,x^pyx^p)^n \ovy\\
\cdot(\ovx^{\,p}\,\ovy\,\ovx^{\,p+1}\,\ovy\,\ovx^{\,p}\,\ovy^{\,q+1})^n
\ovx^{\,p}\,\ovy\,\ovx^{\,p} (x^{p}yx^py^{q}x^pyx^p\,\ovy\,
\ovx^{\,p}\,\ovy)^n y;
\end{multline*}

\item[(viii)]
$m_0(p,q)\equiv x^p y x^p\,\ovy\,\ovx^{\,p}\,\ovy\,x^p y x^p
y^q$;

\medskip

\item[(ix)]
$m_1(p,q)\equiv x^p y x^{p+1} y x^p y^{q+1}$;

\medskip

\item[(x)]
$\partial P(p,q)\equiv x^p y x^{2p+1} y x^p y^q x^p y x^p y^q
$.
\hfill\qed
\een
\end{lem}

\begin{lem}\label{hyper2}
The knot $K(p,q)$ is trivial if $(p,q)=(0,0),(0,1)$, and a
torus knot if $p=0$ and $q\neq -1,0$ (with $K(0,q)=T(2,2q-1)$)
or $(p,q)=(-1,1)$ (with $K(-1,1)=T(5,8))$. In all other cases,
$K(p,q)$ is a hyperbolic tunnel number one knot which is not
2-bridge, and its Seifert Klein bottle $P(p,q)$ is
$\pi_1$-injective in the knot exterior iff $(p,q)$ is not a
pair of the form $(-1,2), \ (-2,1)$, or $(p,0)$.
\end{lem}

\begin{figure}
\begin{center}
\psfrag{H}{$H$}\psfrag{x}{$x$}\psfrag{y}{$y$}
\Figw{4in}{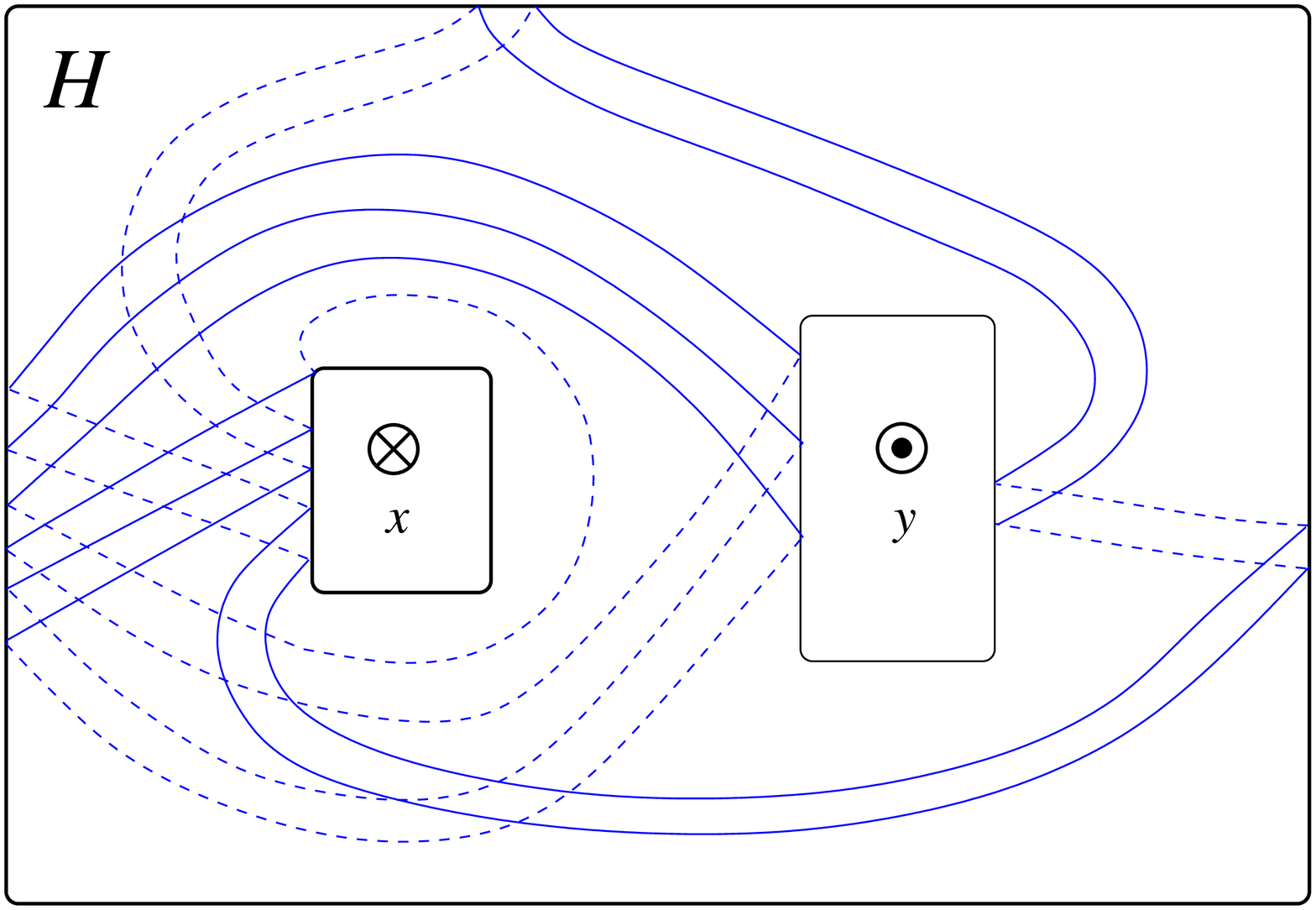}{The knot
$K(-1,1)\subset\partial H$.}{a11b2}
\end{center}
\end{figure}

\bpf
It is easy to see that $K(0,0)$ and $K(0,1)$ are trivial knots,
$K(-1,1)=T(5,8)$ (see Fig.~\ref{a11b2}), and $K(0,q)=T(2,2q-1)$
for all $q$. From now on we assume that $(p,q)$ is not a pair
of the form $(0,0),(0,1),(-1,1)$, or $(0,q)$.

By Lemma~\ref{tunnel}, since $\ell'_0(0,p,q)\equiv x(x^py)^2$
is primitive in $\pi_1(H')$, all the knots $K(p,q)$ have tunnel
number one.

By Corollary~\ref{pi1}, the Seifert Klein bottle $P(p,q)$ is
$\pi_1$-injective iff $[\partial P(p,q)]$ is not primitive nor
a power in $\pi_1(H')$. By Lemma~\ref{words}(x), $\partial
P(p,q)\equiv x^p y x^{2p+1} y x^p y^q x^p y x^p y^q, $ so
$\partial P(p,0)\equiv x(x^{2p}y)^3$ is a primitive word.
Suppose now $q\neq 0$. If $p\neq -1$ then by Lemma~\ref{cohen}
$\partial P(p,q)$ is primitive/power iff $q=1$; since $\partial
P(p,1)\equiv x^{p+1}(x^py)^5$, it follows that $\partial
P(p,1)$ is primitive/power iff $|p+1|\leq 1$ iff $p=-2$ (as
$p\neq 0,-1$), in which case $\partial P(-2,1)$ is primitive.
If $p=-1$ then $\partial
P(-1,q)=\ovx\,y\,\ovx\,y\,\ovx\,y^q\,\ovx\,y\,\ovx\,y^q$, so by
Lemma~\ref{cohen}, as $q\neq 0$, $\partial P(-1,q)$ is a power
iff $q=1$ (which corresponds to the torus knot
$K(-1,1)=T(5,8)$) and primitive iff $q=2$. Therefore, $P(p,q)$
is $\pi_1$-injective iff $(p,q)$ is not one of the pairs
$(-1,2),\ (-2,1)$, or $(p,0)$.

Now, for $p\neq 0$, the circle $m_1(p,q)$ is a positive or
negative braid on $s=3$ strings with $c=6|p|$ crossings (see
Fig.~\ref{figmer}), whence its genus is $g=(c-s+1)/2=3|p|-1\geq
2$ (cf \cite{stallings}), and so the meridian circle of
$P(p,q)$ is a nontrivial knot which is distinct from the genus
one knot $T(2,3)$. Since $\partial P(-1,2), \ \partial P(-2,1),
\ \partial P(p,0)$ are all primitive in $\pi_1(H')$ by the
argument above, the hyperbolicity of $K(p,q)$ for
$(p,q)\neq(0,0),(0,1),(-1,1),(0,q)$ can be established via
Lemma~\ref{hyper1} by checking that
$\mu(\ell'_0(n,p,q))=1=\mu(\ell'_1(n,p,q))$ and either
$\mu(m_0(p,q))=1$ or $\mu(m_1(p,q))=1$ hold in $H'=X(P(p,q))$
for all $n$ and the allowed values of $(p,q)$.

By Lemma~\ref{words}, $m_0(p,q)\equiv x^p y
x^p\,\ovy\,\ovx^{\,p}\,\ovy\,x^p y x^p y^q\in\pi_1(H')$. As
$p\neq 0$, if $q\neq 0$ then the word for $m_0(p,q)$ is
cyclically reduced, while the cyclic reduction of $m_0(p,0)$ is
$x^{2p} y x^p\,\ovy\,\ovx^{\,p}\,\ovy\,x^p y$; in either case
the word for $m_0(p,q)$ contains both $y$ and $\ovy$ factors
and so it is neither primitive nor a power by
Lemma~\ref{cohen}; thus $\mu(m_0(p,q))=1$ for all $p,q$ by
Lemma~\ref{comp}(b).

Consider now the word $$\ell'_0(n,p,q)\equiv
(x^{p+1}yx^py^{q+1}x^py)^n x^{p+1}y (\ovy\,
x^{p}yx^py^{q}x^pyx^p\,\ovy\, \ovx^{\,p})^n x^py\in\pi_1(H').$$
Clearly $\ell'_0(0,p,q)\equiv x(x^py)^2$ is primitive. If
$|n|\geq 2$ then, since $p\neq 0$, it is not hard to see that
the cyclic reduction of the word $\ell'_0(n,p,q)$ contains both
$y$ and $\ovy$ factors, while for $n=\pm 1$ the word
$\ell'_0(n,p,q)$ contains both $y^n$ and $y^{3n}$ factors.
Hence in all cases $\mu(\ell'_0(n,p,q))=1$ by
Lemmas~\ref{cohen} and \ref{comp}(b).

Finally,
$\ell'_1(0,p,q)\equiv\ovx^{\,p}\,\ovy\,\ovx^{\,p}\,\ovy^{\,q}\,
\ovx^{\,p}\,\ovy\,\ovx^{\,p}\,  \ovy^{\,q-1}\in\pi_1(H')$.
For $n=0$, if $|p|\geq 2$ then $\ell'_1(0,p,q)$ contains at
least three factors $\ovx\,^p$; if $|p|=1$ then
$\ell'_1(0,p,0)$ contains both $y$ and $\ovy$ factors,
$\ell'_1(0,p,1)=\ovx\,^p(\ovx\,^p\,y)^3$ is primitive, and for
any $q\neq 0,1$ the word $\ell'_1(0,p,q)$ contains three
different powers $\ovy,\ovy\,^q,\ovy\,^{q-1}$. Finally, for
$|n|=1$ we have $\ell'_1(n,p,q)\equiv(x(x^p y)^2)^{|n|}$, a
primitive word, and for $|n|\geq 2$ the word $\ell'_1(n,p,q)$
contains both $y$ and $\ovy$ factors. Thus again
$\mu(\ell'_1(n,p,q))=1$ in all cases.
\epf

\begin{lem}\label{toro}
\ben
\ita
The boundary slope of the surface $P(p,q)$ is $r=4q-36p$.
\itb
If $(p,q)$ is not a pair of the form $(0,0),(0,1),(-1,1),(0,q)$
then, except for $K(-1,2)(r)=S^2(2,2,3)$,
$K(-2,1)(r)=S^2(2,2,7)$, and $K(p,0)(r)=S^2(2,2,|6p-1|)$, the
manifold\\ $K(p,q)(r)$ is irreducible and toroidal.
\een
\end{lem}

\bpf
Since $\partial P(p,q)$ lies in $\partial H$, the boundary
slope of  $\partial P(p,q)$ coincides with the linking number
between $\partial P(p,q)$ and a parallel copy in $\partial H$.
From Fig.~\ref{fig11b} we can thus see that the boundary slope
of $\partial P(p,q)$ is $r=0$, and that the above linking
number decreases by $6^2=36$ with each positive Dehn twist
along $D_0$, and increases by $2^2=4$ with each positive Dehn
twist along $D_1$. Thus the boundary slope of $P(p,q)$ is
$r=4q-36p$ and (a) holds.

For part (b) suppose $(p,q)\neq(0,0),(0,1),(-1,1),(0,q)$, and
denote $P(p,q)$ by $P$ for simplicity. Let $\wt{X}$ and
$\wt{N}$ be the manifolds obtained by attaching 2-handles to
$H'=X(P)$ and $N(P)$ along $\partial P$, respectively; thus, if
$r$ is the boundary slope of $P$ then
$K(p,q)(r)=\wt{X}\cup_{\partial}\wt{N}$. We will consider
$N(P)$ as a Seifert fibered space over a disk with two singular
fibers of indices 2,2; in particular, the pair
$(\wt{N},\partial\wt{N})$ is irreducible, and the circle
$\ell'_0(0,p,q)\subset \partial\wt{N}$ is a fiber of $\wt{N}$
disjoint from $\partial P(p,q)$.

By Lemma~\ref{hyper2}, if $(p,q)\neq(-1,2),(-2,1),(p,0)$ then
the surface $P$ is $\pi_1$-injective in the exterior of
$K(p,q)$, which implies that the surface $T_{P}$ is
incompressible in $H'=X(P)$ and hence, by the 2-handle addition
theorem, that the pair $(\wt{X},\partial\wt{X})$ is
irreducible. Therefore $K(r)=\wt{X}\cup_{\partial}\wt{N}$ is
irreducible and $\partial\wt{X}=\partial\wt{N}$ is an
incompressible torus in $K(r)$.

Now, by Lemma~\ref{words}, $[\partial P(p,0)]\equiv
x(x^{2p}y)^3$ is primitive in $\pi_1(H')=\pi_1(X(P))$ and so
$\wt{X}$ is a solid torus; thus the Seifert fibration of
$\wt{N}$ extends to a Seifert fibration of
$K(p,0)(r)=\wt{X}\cup_{\partial}\wt{N}$ over the 2-sphere with
fibers of indices $2,2,n$, where $n=|6p-1|$ is the order of the
cyclic group
$$
\pi_1(\wt{X})/\grp{[\ell'_0(0,p,0)]}=
\pi_1(H')/\grp{[\partial
P(p,0)],[\ell'_0(0,p,0)]}=\mathbb{Z}_{|6p-1|}.
$$
Therefore $K(p,0)(r)=S^2(2,2,|6p-1|)$. The identities
$K(-1,2)(r)=S^2(2,2,3)$ and $K(-2,1)(r)=S^2(2,2,7)$ follow in a
similar way.
\epf

We now classify the words $[c_i'(n,p,q)]$ and $[c_i''(n,p,q)]$
which are primitive or a power in $\pi_1(H')$; we will need
this information to discern which knots $K(p,q)$ admit a (1,1)
decomposition.

\begin{lem}\label{wc}
In $\pi_1(H')$, for $p\neq 0$,
\ben
\ita
$\bullet$ $c_0'(n,p,q)$ is primitive iff $(n,p,q)=(0,\pm 1,q),$
$(1,-2,0)$, $(1,-1,1)$,  $(1,p,-1)$, $(-1,2,0)$, $(-1,1,1)$,
$(2,-1,0)$, $(-2,1,0)$, or $(n,p,q)= (n,1,-1)$ for all $n$;

$\bullet$ $c_0'(n,p,q)$ is a power iff $(n,p,q)=(0,p,q)$ for
$|p|\geq 2$, or $(n,p,q)=(1,-1,0)$, $(-1,1,0)$;

\medskip

\itb
$\bullet$ $c_0''(n,p,q)$ is primitive iff $(n,p,q)=(0,2,q),$
$(1,p,0)$, $(1,-1,2)$, $(1,-2,1)$;

$\bullet$ $c_0''(n,p,q)$ is a power iff $(n,p,q)=(0,1,q)$,
$(1,-1,1)$;

\medskip

\itc
$\bullet$ $c_1'(n,p,q)$ is primitive iff $(n,p,q)=(0,\pm 1,1)$,
$(0,\pm 1,3)$, $(-1,2,0)$, $(-1,-1,2)$, or $(-1,p,1)$;

$\bullet$ $c_1'(n,p,q)$ is a power iff $(n,p,q)=(-1,1,0)$ or
$(0,p,2)$;

\medskip

\itd
$\bullet$ $c_1''(n,p,q)$ is primitive iff $(n,p,q)=(0,\pm 1,2)$
or $(0,p,0)$;

$\bullet$ $c_1''(n,p,q)$ is a power iff $(n,p,q)=(0,p,1)$.
\een
\end{lem}

\bpf
The words for $c'_i(n,p,q)$ and $c''_i(n,p,q)$ in
$\pi_1(H')=\grp{x,y \ | \ -}$ are given in Lemma~\ref{words}.
It is easy to see that in the given cases the words
$c'_i(n,p,q)$ and $c''_i(n,p,q)$ are indeed primitive or powers
as claimed. In order to establish the converse statements we
apply Lemma~\ref{cohen}; the fact that a word of the form
$u^sv^t$ in the free group $\grp{u,v \ | \ -}$ is primitive iff
$|s|=1$ or $|t|=1$, and a power iff $s=0,|t|\geq 2$ or $|s|\geq
2,t=0$, will also be of use.

So suppose some word $c'_i(n,p,q)$ or $c''_i(n,p,q)$ is
primitive or a power in $\pi_1(H')$. We consider the case of
the word $c_0'(n,p,q)\equiv x^p(x^{p+1}yx^py^{q+1}x^py)^n$ in
part (a) in full detail; the other parts of the lemma follow
along entirely similar lines, so their proof will be omitted.

\medskip

(1) \ For $n=0$, $c_0'(0,p,q)\equiv x^p$. Thus $c_0'(0,p,q)$ is
primitive iff $p=\pm 1$ and a power iff $|p|\geq 2$.

\medskip

(2) \ For $n=1$, $c_0'(1,p,q)\equiv x^{2p+1}yx^py^{q+1}x^py$.
\bit
\item
If $q\neq -1$ then $x^{2p+1}yx^py^{q+1}x^py$ is cyclically
reduced, so by Lemma~\ref{cohen} either all exponents of $x$
are $1$ or all $-1$, or all exponents of $y$ are $1$ or all
$-1$; that is, either $p=-1$ or $q=0$. By Lemma~\ref{cohen},
the word $c_0'(1,-1,q)\equiv
\ovx\,y\ovx\,y^{q+1}\,\ovx\,y$ is
primitive iff $q=1$ and a power iff $q=0$, while the word
$c_0'(1,p,0)\equiv x^{2p+1}yx^pyx^py=x^{p+1}(x^py)^3$ is
primitive iff $p+1=\pm 1$ iff $p=-2$ (as $p\neq 0$) and a power
iff $p=-1$.

\item
If $q=-1$ then $c_0'(1,p,-1)\equiv x^{2p+1}yx^{2p}y$ is
primitive for all $p$.
\eit

{\bf Remark:} for the rest of the argument in this proof, we
will implement the strategy used in (2) and will not explicitly
indicate the use of Lemma~\ref{cohen} for the sake of brevity.

(3) \ For $n=-1$,  $c_0'(-1,p,q)\equiv xyx^py^{q+1}x^py$.
\bit
\item
If $q\neq -1$ then $xyx^py^{q+1}x^py$ is cyclically reduced,
hence either $p=1$ or $q=0$. The word $c_0'(-1,1,q)\equiv
xyxy^{q+1}xy$ is primitive iff $q=1$ and a power iff $q=0$,
while the word $c_0'(-1,p,0)\equiv xyx^pyx^py$ is primitive iff
$p=2$ and a power iff $p=1$.

\item
If $q=-1$ then the word $c_0'(-1,p,-1)\equiv xyx^{2p}y\equiv
x^{2p-1}(xy)^2$ is primitive iff $p=1$ (as $p\neq 0$) and never
a power.

\eit

(4) \  For $n\geq 2$, $c_0'(n,p,q)\equiv
x^p(x^{p+1}yx^py^{q+1}x^py)^n$.
\bit
\item
If $p\neq -1$ and $q\neq -1$ then $c_0'(n,p,q)\equiv
x^p(x^{p+1}yx^py^{q+1}x^py)^n$ is cyclically reduced and $x$
appears with the three distinct exponents $2p+1,p+1,p$, so
$c_0'(n,p,q)$ is never primitive nor a power.

\item
If $p=-1$, $c_0'(n,-1,q)\equiv
\ovx(y\ovx\,y^{q+1}\,\ovx\,y)^n$ is cyclically reduced and $y$
appears with exponents $1,2,q+1$, hence we must have
$q+1=0,1,2$, ie $q=-1,0,1$.
\bit
\item[--]
If $q=-1$ then $c_0'(n,-1,-1)\equiv
\ovx(y\,\ovx\,^{2}\,y)^n=\ovx\,y\,\ovx\,^{2}\,y^2\,
\ovx\,^{2}\,y(y\,\ovx\,^{2}\,y)^{n-2}$ is cyclically reduced
and both $\ovx$ and $y$ appear with exponents 1 and 2, so
$c_0'(n,-1,-1)$ is never primitive nor a power.
\item[--]
If $q=0$ then $c_0'(n,-1,0)\equiv
\ovx(y\ovx\,y\,\ovx\,y)^n\equiv u^{n-1}v^3$ for the basis
$u=y\ovx\,y\,\ovx\,y$ and $v=y\ovx$ of $\pi_1(H')$, hence
$c_0'(n,-1,0)$ is primitive iff $n=2$ and never a power.
\item[--]
If $q=1$ then $c_0'(n,-1,1)\equiv
\ovx(y\ovx\,y^{2}\,\ovx\,y)^n\equiv u^{2n-1}v^2$ for the basis
$u=y\ovx\,y$ and $v=y\ovx$ of $\pi_1(H')$, hence $c_0'(n,-1,1)$
is neither primitive nor a power.
\eit

\item
For $q=-1$ and $p\neq -1$ the word $$c_0'(n,p,-1)\equiv
x^p(x^{p+1}yx^{2p}y)^n = x^{2p+1}yx^{2p}y x^{p+1}yx^{2p}y
(x^{p+1}yx^{2p}y)^{n-2}$$ is cyclically reduced and $x$ appears
with exponents $2p+1,2p,p+1$, which are mutually distinct for
$p\neq 1$, in which case $c_0'(n,p,-1)$ is neither primitive
nor a power, while for $p=1$ the word $c_0'(n,1,-1)\equiv
x(x^{2}y)^{2n}$ is always primitive.
\eit

(5) \  For $n\leq -2$,
\begin{align*}
c_0'(n,p,q) & \equiv
\ovx\,^p(x^{p+1}yx^py^{q+1}x^py)^{|n|}\\
& =
xyx^py^{q+1}x^py x^{p+1}yx^py^{q+1}x^py
(x^{p+1}yx^py^{q+1}x^py)^{|n|-2}.
\end{align*}
\bit
\item
If $q\neq -1$ then the word for $c_0'(n,p,q)$ is cyclically
reduced and $x$ appears with exponents $p+1,p,1$ while $y$
appears with exponents $q+1,1$, so we must have $p=1$ and
$q=0$; thus $c_0'(n,1,0)\equiv\ovx(x^{2}yxyxy)^{|n|}\equiv
u^{|n|-1}v^3$ for the basis $u=x^{2}yxyxy$ and $v=xy$ of
$\pi_1(H')$, so $c_0'(n,1,0)$ is primitive iff $n=-2$ and never
a power.

\item
If $q=-1$ then $c_0'(n,p,q) \equiv xyx^{2p}y x^{p+1}yx^{2p}y
(x^{p+1}yx^{2p}y)^{|n|-2}$ and $x$ appears with exponents
$2p,p+1,1$, so $p=1$ again; clearly
$c_0'(n,1,-1)\equiv\ovx(x^{2}y)^{2|n|}$ is always primitive.
\eit

Therefore part (a) holds.
\epf

The following lemma is the last result needed in the proof of
Theorem~\ref{main}.

\begin{lem}\label{one}
If $K(p,q)$ is a hyperbolic knot then $K(p,q)$ is a
$(1,1)$-knot iff $(p,q)$ is a pair of the form $(-1,0)$,
$(1,q)$, or $(p,1)$, $(p,2)$ for $p\neq 0$.
\end{lem}

\bpf
By Lemma~\ref{hyper2}, $(p,q)$ is not a pair of the form
$(0,0)$, $(0,1)$, $(-1,1)$, or $(0,q)$. That $K(p,q)$ admits a
(1,1) decomposition for $(p,q)=(p,1)$, $(p,2)$, $(1,q)$, and
$(-1,0)$ follows directly from Lemma~\ref{lem3}(a) using the
following choices for $\wh{c}_0,\wh{c}_1$ and the fact that
$H'=X(P(p,q))$:

\ben
\item
$\wh{c}_0=c'_0(0,p,1)$ and $\wh{c}_1=c''_1(0,p,1)$; then
$[\wh{c}_0]\equiv x^p$, $[\wh{c}_1]\equiv (yx^p)^2$, whence
$\mu(\wh{c}_0)=|p|$ and $\mu(\wh{c}_1)=2$, and $[\partial
D_2(0,p,1)]\equiv(yx^p)^2(\ovy\,\ovx\,^p)^2\equiv[x^{|p|},v^2]$
for the basis $\{x,v=yx^p\}$ of $\pi_1(H')$;

\item
$\wh{c}_0=c'_0(0,p,2)$ and $\wh{c}_1=c'_1(0,p,2)$; then
$[\wh{c}_0]\equiv x^p$, $[\wh{c}_1]\equiv (yx^p)^2$, so again
$\mu(\wh{c}_0)=|p|$ and $\mu(\wh{c}_1)=2$ while $[\partial
D_2(0,p,2)]\equiv(yx^p)^2(\ovy\,\ovx\,^p)^2\equiv[x^{|p|},v^2]$
for the basis $\{x,v=yx^p\}$ of $\pi_1(H')$;

\item
$\wh{c}_0=c''_0(0,1,q)$ for $q\neq 1,2$ (the cases $q=1,2$
follow from (1) and (2) above); then $[c''_0(0,1,q)]=(yx)^2$,
so $\mu(\wh{c}_0)=2$, while $\mu(\wh{c}_1)=1$ by
Lemma~\ref{wc}(c) and $[\partial D_2(0,1,q)]\equiv [u^2,y]$ for
the basis $\{u=yx,y\}$ of $\pi_1(H')$;

\item
$\wh{c}_0=c'_0(1,-1,0)$; then $[\wh{c}_0]\equiv (y\,\ovx)^3$,
so $\mu(\wh{c}_0)=3$, and $\mu(\wh{c}_1)=1$ by
Lemma~\ref{wc}(c) while $[\partial D_2(1,-1,0)] \equiv
(y\,\ovx)^2y(x\,\ovy)\,^3\ovx\equiv[u^3,y]$ for the basis
$\{u=y\,\ovx,y\}$ of $\pi_1(H')$.
\een

Conversely, assume that $K(p,q)$ admits a (1,1) decomposition
and the pair $(p,q)$ is not of the form $(-1,0)$, $(1,q)$, or
$(p,1)$, $(p,2)$ for $p\neq 0$ (recall also that $(p,q)$ is not
a pair of the form $(0,0)$, $(0,1)$, $(-1,1)$, or $(0,q)$); we
show this situation contradicts Lemma~\ref{lem3}(a). The
following fact will be useful in the sequel.

\begin{claim}\label{n}
If $\mu(\wh{c}_0)=1=\mu(\wh{c}_1)$ then $n\neq -1,0$.
\end{claim}
For $n=-1,0$ the word for $[\partial D_2(n,p,q)]$ in
$\pi_1(H')$ has the following cyclic reductions:
\begin{align*}
[\partial D_2(0,p,q)]\rule{8pt}{0pt}  & \equiv(yx^p)^2(\ovy\,\ovx\,^p)^2; \\
[\partial D_2(-1,p,q)]  & \equiv
x^pyx^py^qx^pyx^p\ovy\,\ovx\,^p\,\ovy\,\ovx\,\ovy\,\ovx\,^p\,\ovy\,^q\,\ovx\,^p\,
\ovy\,\ovx\,^p\,\ovy\,^q\,\ovx\,^p\,\ovy\,\ovx\,^p\,yx^py^qyxyx^py^q
\text{ for } q\neq 0;\\
[\partial D_2(-1,p,0)]  & \equiv
x^{2p}yx^{2p}yx^p\,\ovy\,\ovx\,^p\,\ovy\,\,\ovx\,
\ovy\,\ovx\,^{2p}\,\ovy\,\ovx\,^{2p}\,\ovy\,\ovx\,^p\,yx^pyxy.
\end{align*}
Thus $[\partial D_2(n,p,q)]\not\equiv[x,y]$ for $n=-1$ or
$n=0$, so the claim follows by Lemma~\ref{lem3}(a)(i).
\hfill\qed(Claim~\ref{n})

We now consider two cases:

\begin{case}
$\mu(\wh{c}_0)=1$.
\end{case}

By Lemmas~\ref{comp}(b) and \ref{lem3}(b), in $\pi_1(H')$,
$[\wh{c}_0]$ is not a power and $[\wh{c}_1]$ is either
primitive with $\mu(\wh{c}_1)=1$ or a power with
$\mu(\wh{c}_1)\geq 2$. We thus consider the following subcases.

\begin{subcase}
$[\wh{c}_1]$ is primitive, $\mu(\wh{c}_1)=1$.
\end{subcase}

As $\wh{c}_1=c_1'(n,p,q)$ or $c_1''(n,p,q)$ gives rise to a
primitive word in $\pi_1(H')$ then $n=-1,0$ by Lemma~\ref{wc},
contradicting Claim~\ref{n}.

\begin{subcase}
$[\wh{c}_1]$ is a power, $\mu(\wh{c}_1)\geq 2$.
\end{subcase}

By Lemma~\ref{wc}, $\wh{c}_1$ must be one of $c'_1(-1,1,0)$,
$c'_1(0,p,2)$, or $c''_1(0,p,1)$, contradicting our hypothesis
that $(p,q)$ is not of the form $(1,q)$, $(p,1)$, $(p,2)$.

\begin{case}
$\mu(\wh{c}_0)\geq 2$.
\end{case}

Then $[\wh{c}_0]$ is a power in $\pi_1(H')$; by Lemma~\ref{wc},
with the given restrictions on the pair $(p,q)$, we must have
$\wh{c}_0=c'_0(0,p,q)$ for $|p|\geq 2$ and $q\neq 1,2$, and
hence either $[\wh{c}_1]=[c'_1(0,p,q)]\equiv y^{q-1}x^pyx^p$ or
$[\wh{c}_1]=[c''_1(0,p,q)]\equiv y^{q}x^pyx^p$.

Now, by Lemma~\ref{cohen}, for any $p\neq 0$, $y^{q-1}x^pyx^p$
is a power iff $q=2$ while $y^{q}x^pyx^p$ is a power iff $q=1$;
thus $\mu(\wh{c}_1)=1$ in all cases. However,
$\mu(\wh{c}_0)=|p|\geq 2$ since $[\wh{c}_0]=[c'_0(0,p,q)]\equiv
x^p$, while $[\partial D_2(0,p,q)]\equiv
(yx^p)^2(\ovy\,\ovx\,^p)^2\not\equiv[x^p,y]$ for the basis
$\{x,y\}$ of $\pi_1(H')$, contradicting
Lemma~\ref{lem3}(a)(ii).
\epf



\end{document}